\documentclass[11pt]{amsart}
\usepackage{amsmath,amssymb,amscd,amsthm}
\newtheorem{formula}{}[section]
\newtheorem{proposition}[formula]{Proposition}
\newtheorem{corollary}[formula]{Corollary}
\newtheorem{lemma}[formula]{Lemma}
\newtheorem{theorem}[formula]{Theorem}
\theoremstyle{definition}
\newtheorem{definition}[formula]{Definition}
\newtheorem{example}[formula]{Example}
\theoremstyle{remark}
\newtheorem*{remark}{Remark}

\newcommand{\D}{\Delta}
\newcommand{\e}{\varepsilon}
\newcommand{\C}{\mathbb C}
\newcommand{\R}{\mathbb R}
\newcommand{\Q}{\mathbb Q}
\newcommand{\Z}{\mathbb Z}
\renewcommand{\k}{\mathbf k}

\newcommand{\bideg}{\mathop{\rm bideg}}
\newcommand{\bs}{\mathop{\rm bs}}
\newcommand{\cc}{\mathop{\rm cc}}
\newcommand{\cub}{\mathop{\rm cub}}
\newcommand{\Tor}{\mathop{\rm Tor}\nolimits}
\newcommand{\cone}{\mathop{\rm cone}}
\newcommand{\ma}{\mathop{\rm ma}}
\newcommand{\zk}{\mathcal Z_K}
\newcommand{\wk}{\mathcal W_K}

\newcommand{\ep}{\nolinebreak\hfill$\square$}

\begin{document}

\title[Moment-angle complexes and simplicial manifolds]
{Moment-angle complexes and combinatorics of simplicial manifolds}
\author{Victor M. Buchstaber}
\author{Taras E. Panov}
\thanks{Partially supported by the Russian Foundation for Fundamental
Research, grant no. 99-01-00090}
\address{Department of Mathematics and Mechanics, Moscow
State University, 119899 Moscow, RUSSIA}
\email{buchstab@mech.math.msu.su \quad tpanov@mech.math.msu.su}

\begin{abstract}
Let $\rho:(D^2)^m\to I^m$ be the orbit map for the diagonal action of the
torus $T^m$ on the unit poly-disk $(D^2)^m\subset\C^m$, $I^m=[0,1]^m$ is the
unit cube. Let $C$ be a cubical subcomplex in $I^m$. The moment-angle complex
$\ma(C)$ is a $T^m$-invariant bigraded cellular decomposition of the subset
$\rho^{-1}(C)\subset(D^2)^m$ with cells corresponding to the faces of $C$.
Different combinatorial problems concerning cubical complexes and related
combinatorial objects can be treated by studying the equivariant topology of
corresponding moment-angle complexes. Here we consider moment-angle complexes
defined by canonical cubical subdivisions of simplicial complexes. We
describe relations between the combinatorics of simplicial complexes and the
bigraded cohomology of corresponding moment-angle complexes. In the case when
the simplicial complex is a simplicial manifold the corresponding moment-angle
complex has an orbit consisting of singular points. The complement of an
invariant neighbourhood of this orbit is a manifold with boundary. The
relative Poincar\'e duality for this manifold implies the generalized
Dehn--Sommerville equations for the number of faces of simplicial manifolds.
\end{abstract}

\maketitle

\section{Introduction}

The classical Dehn--Sommerville equations for simplicial convex
$n$-dim\-en\-si\-onal polytope $P^n$ give a set of linear relations among the
numbers $f_i$ of $i$-dimensional faces of $P^n$. The integer vector
$f(P^n)=(f_0,f_1,\ldots,f_{n-1})$ is called the $f$-{\it vector} of $P^n$. We
put $f_{-1}=1$. The Dehn--Sommerville equations were established by Dehn for
$n\le5$, and by Sommerville in the general case in 1927 (see~\cite{So}).
Their original form is as follows
\begin{equation}
\label{DSf}
  f_k=\sum_{j=k}^{n-1}(-1)^{n-1-j}\binom{j+1}{k+1}f_j,\quad k=0,\ldots,n-1.
\end{equation}
Define the $h$-{\it vector} from the equation
\begin{equation}
\label{hvector}
  h_0t^n+\ldots+h_{n-1}t+h_n=(t-1)^n+f_0(t-1)^{n-1}+\ldots+f_{n-1}.
\end{equation}
Obviously, the $f$-vector and the $h$-vector determine each other by means of
linear equations (note that $h_0=1$). For instance,
\begin{equation}
\label{hf}
  h_k=\sum_{i=0}^k(-1)^{k-i}\binom{n-i}{k-i}f_{i-1}.
\end{equation}
The notion of $h$-vector gives rise to the simplest and the most elegant
form of the Dehn--Sommerville equations~(\ref{DSf}) for simplicial polytopes:
\begin{equation}
\label{DSpol}
  h_i=h_{n-i},\quad i=0,\ldots,n.
\end{equation}

The Dehn--Sommerville equations can be generalized quite widely.
In~\cite{Kl} V.~Klee reproved the Dehn--Sommerville equations in the
form~(\ref{DSf}) in a more general context of {\it Eulerian manifolds}. In
particular, it turns out that equations~(\ref{DSf}) hold for any simplicial
manifold (i.e. triangulated topological manifold) of dimension $n-1$.
Analogues of equations~(\ref{DSf}) were obtained by Bayer and
Billera~\cite{BB} (for Eulerian partially ordered sets) and Chen and
Yan~\cite{CY} (for arbitrary polyhedra).

It follows directly from~(\ref{DSpol}) that the affine hull of $f$-vectors
$(f_0,\ldots,f_{n-1})$ of simplicial polytopes in $n$-dimensional space is an
$\bigl[\frac n2\bigr]$-dimensional plane. The same is true for the affine
hull of vectors $(b_0,b_1,\ldots,b_n)$ of Betti numbers $b_i=\dim H_i(M)$ of
orientable connected (simplicial) $n$-dimensional manifolds $M^n$ due to the
Poincar\'e duality:
\begin{equation} \label{pd}
  b_i(M^n)=b_{n-i}(M^n),\quad i=0,\ldots,n.
\end{equation}
This similarity between $f$-vectors of simplicial polytopes and Betti numbers
of orientable manifolds was pointed out by Klee in~\cite{Kl}. Moreover, it
turns out that the ``combinatorial" duality~(\ref{DSpol}) can be
interpreted in terms of the ``topological" duality~(\ref{pd}). Given an
$n$-dimensional simplicial polytope $P$ (or, more generally, a complete
simplicial fan $\Sigma$), one may construct the {\it toric variety} $M_P$ (or
$M_\Sigma$) of (real) dimension $2n$ (see~\cite{Da},~\cite{Fu}). This variety
is not necessarily a manifold, however its homology (with rational
coefficients in general) satisfies the Poincar\'e duality, and its even Betti
numbers equal the components of $h$-vector of $P$: $b_{2i}(M_P)=h_i(P)$. This
gives a ``topological" proof of the Dehn--Sommerville
equations~(\ref{DSpol}).

The dual (or polar) to any simplicial polytope $P$ is a {\it simple}
polytope, which we denote $P^{*}$. Given an $n$-dimensional simple polytope
$P^{*}$ with $m$ facets, Davis and Januszkiewicz defined in~\cite{DJ} a
manifold $\mathcal Z_{P^{*}}$ of dimension $m+n$ acted on by the torus $T^m$.
This manifold depends only on the combinatorial type (i.e., the face lattice)
of a polytope, and the orbit space for the $T^m$-action is combinatorially
$P^{*}$. The manifolds $\mathcal Z_{P^{*}}$ establish the bridge between
topology of manifolds and combinatorics of polytopes (or more general
objects, such as simplicial spheres). Various connections between
topology of $\mathcal Z_{P^{*}}$ and toric geometry, symplectic geometry,
subspace arrangements, and combinatorial theory of $f$-vectors were studied
in~\cite{BP1}, \cite{BP2}, \cite{BP3}. Any smooth toric variety
(or Hamiltonian $T^n$-manifold) defined by a simple polytope of combinatorial
type $P^{*}$ is the quotient of $\mathcal Z_{P^{*}}$ by a freely acting toric
subgroup $T^{m-n}\subset T^m$. This is also true for topological analogues of
smooth toric varieties, which we call {\it quasitoric manifolds}, also
introduced in~\cite{DJ}. For more information about quasitoric manifolds and
their topology see~\cite{BaBe}, \cite{BR1}, \cite{BR2}, \cite{Ma},
\cite{Pa1}, \cite{Pa2}. We showed in~\cite{BP2}, \cite{BP3} that the
cohomology algebra of $\mathcal Z_{P^*}$ with coefficients in any field $\k$
is isomorphic to the $\Tor$-algebra $\Tor_{\k[v_1,\ldots,v_m]}\bigl(\k(P),\k
\bigr)$, where $P$ is polar to $P^{*}$, $m=f_0$ is the number of vertices of
$P$, and $\k(P)$ is the {\it
Stanley--Reisner face ring} of $P$. The Koszul resolution gives then a very
simple model $H\bigl[\k(P)\otimes\Lambda[u_1,\ldots,u_m],d\bigr]$ for the
cohomology algebra of $\mathcal Z_{P^*}$. In particular, the cohomology of
$\mathcal Z_{P^*}$ acquires a {\it bigraded} algebra structure, and the
bigraded Poincar\'e duality implies the Dehn--Sommerville
equations~(\ref{DSpol}).

The boundary complex of a simplicial polytope is a simplicial sphere.
However, now it is well known that not any triangulation of a topological
sphere can be obtained in such way. First examples of such ``non-polytopal"
spheres were found by Gr\"unbaum, and the smallest non-polytopal sphere is of
dimension 3 with 8 vertices (the so-called {\it Barnette sphere},
see~\cite{Ba}). However, the Dehn--Sommerville relations in the
form~(\ref{DSpol}) still hold for any simplicial sphere
(see~\cite[\S~II.6]{St2}). In~\cite{BP4} we extend our approach to manifolds
$\mathcal Z_{P^{*}}$ defined by simple polytopes to the case of arbitrary
simplicial complex $K$. We endow the cone $\cone(K)$ with a structure of
cubical complex and embed it into the boundary complex of $m$-dimensional
unit cube $I^m$ (here $m=f_0(K)$ is the number of vertices of $K$). Then we
view $I^m$ as the orbit space for the diagonal action of the torus
$T^m=S^1\times\cdots\times S^1$ on the unit poly-disk $(D^2)^m$ in
$m$-dimensional complex space $\C^m$. Hence, the above cubical embedding
$\cone(K)\to I^m$ is covered by a $T^m$-equivariant embedding
$\zk\to(D^2)^m$, where $\zk$ is a cellular complex canonically decomposed
into the union of blocks $(D^2)^n\times T^{m-n}$ with $n=\dim K+1$. We call
this $\zk$ the {\it moment-angle complex associated to the simplicial complex
$K$}. If $K$ is a simplicial sphere, then $\zk$ is a manifold. (If, moreover,
$K$ is polytopal, i.e. $K=\partial P$, then $\zk$ coincides with the above
described manifold $\mathcal Z_{P^{*}}$, where $P^{*}$ is polar to $P$.)
However, in general, $\zk$ has more complicated structure. In~\cite{BP4} we
showed that our moment-angle complex $\zk$ is homotopy equivalent to the
complement of complex coordinate subspace arrangement defined by $K$, and its
cohomology algebra is isomorphic, as in the polytopal case, to the
$\Tor$-algebra $\Tor_{\k[v_1,\ldots,v_m]}\bigl(\k(K),\k \bigr)$.  We note
that the Betti numbers of the complement of a {\it real\/} coordinate
subspace arrangement were calculated in terms of resolution of the
Stanley--Reisner ring in~\cite{GPW}. In the case when $K$ is a simplicial
sphere, the bigraded Poincar\'e duality in the cohomology algebra of $\zk$
gives a ``topological proof" of the Dehn--Sommerville equations~(\ref{DSpol})
for simplicial spheres. In this paper we construct a cellular chain complex
that calculates the homology of $\zk$ and gives a very transparent
characterization of homology classes of $\zk$ (and of the complement of a
coordinate subspace arrangement). This chain complex is dual to a certain
cochain subcomplex of the Koszul complex
$\bigl[\k(K)\otimes\Lambda[u_1,\ldots,u_m],d\bigr]$ for the $\Tor$-algebra
$\Tor_{\k[v_1,\ldots,v_m]}\bigl(\k(K),\k \bigr)$.

From the topological viewpoint it is very interesting to study the case when
$K$ is a {\it simplicial manifold}.  The moment-angle complex $\zk$ here is
not a manifold, however, its singularities are tractable. It turns out that
$\zk$ contains the product $|\cone(K)|\times T^m$ of (polyhedron
corresponding to) the cone over $K$ and torus $T^m$. The closure
$W_K=\overline{\zk\setminus|\cone(K)|\times T^m}$ of the complement of this
singular part is a manifold with boundary $|K|\times T^m$. This $W_K$ is
homotopically equivalent to another moment-angle complex $\wk$ which covers
equivariantly the restriction to $K\subset\cone(K)$ of the cubical embedding
$\cone(K)\to I^m$. We construct an appropriate cellular decomposition of
$\wk$, which allows to calculate the homology efficiently. The Poincar\'e
duality for manifolds with boundary gives in this case
$$
  H_i(W_K)\cong H^{m+n-i}(W_K,|K|\times T^m),\quad i=0,\ldots,m+n.
$$
This duality again regards the bigraded structure. As a consequence, we
obtain a topological proof of the Dehn--Sommerville equations~(\ref{DSf})
for simplicial manifolds in the following nice form:
\begin{equation}
\label{DSman}
  h_{n-i}-h_i=(-1)^i\bigl(\chi(K^{n-1})-\chi(S^{n-1})\bigr)\binom ni,
  \quad i=0,1,\ldots,n.
\end{equation}
where $\dim K^{n-1}=n-1$, and $\chi(\cdot)$ denotes the Euler number. We have
$\chi(K^{n-1})=f_0-f_1+\ldots+(-1)^{n-1}f_{n-1}=1+(-1)^{n-1}h_n$ and
$\chi(S^{n-1})=1+(-1)^{n-1}$. This generalizes equations~(\ref{DSpol}) to the
case of arbitrary simplicial manifold. In particular, if $K$ is an
odd-dimensional simplicial manifold, one has $h_{n-i}=h_i$.
In~\cite[(7.11)]{Pa} Pachner proved by means of his bistellar flip theorem
that the value $h_{n-i}-h_i$ is a topological invariant of a PL-manifold
(i.e. is independent on a PL-triangulation). The formula~(\ref{DSman})
calculates this invariant exactly for any simplicial (not necessarily PL)
manifold.

Our moment-angle complexes enable to reformulate many combinatorial
statements and hypotheses concerning $f$-vectors in topological terms. We
have already mentioned the Dehn--Sommerville relations for both simplicial
spheres and simplicial manifolds, however this is only the first and the
simplest example. The most intriguing open problem here is the so-called
$g$-theorem (or McMullen's inequalities) for simplicial spheres~\cite{St3}.
It includes the Generalized Lower Bound hypothesis for simplicial spheres,
which asserts the monotonicity property for the $h$-vector:
\begin{equation}
\label{glb}
  h_0\le h_1\le h_2\le \cdots\le h_{[\frac n2]}.
\end{equation}
For polytopal spheres $g$-theorem was proved by Stanley~\cite{St1}
(necessity) and Billera, Lee~\cite{BL} (sufficiency) in 1980. The first
inequality $h_0\le h_1$ is equivalent to $1\le m-n$, which is obvious. The
second ($h_1\le h_2$, $n\ge4$) is equivalent to the lower bound $f_1\le
nf_0-\binom{n+1}2$ for the number of edges, which is also known for
simplicial spheres. The next inequality $h_2\le h_3$ is still open (for
simplicial spheres). For the history of $g$-theorem and related questions
see~\cite{St2},~\cite{St3},~\cite{Zi}. As we mentioned in~\cite{BP2}, the
inequality $h_1\le h_2$ is equivalent to the upper bound
$b_3(\zk)\le\binom{m-n}2$ for the third Betti number of manifold $\zk$ (we
note that $\zk$ is always 2-connected). Other inequalities from~(\ref{glb})
also acquire such topological interpretation, however, in general case in
terms of {\it bigraded} Betti numbers of $\zk$. We hope that such
inequalities can be deduced from the equivariant topology of~$\zk$.

The authors wish to express special thanks to Oleg Musin for stimulating
discussions and helpful comments, in particular, for drawing our attention to
the results of~\cite{CY} and~\cite{Pa}.

\section{Cubical complexes determined by a simplicial complex}

Let $K^{n-1}$ be an $(n-1)$-dimensional simplicial complex on the vertex
set $[m]=\{1,\ldots,m\}$ (hence, any simplex of $K$ has at most $n$
vertices). As usual, we view $K$ as a set of subsets of $[m]$, that is,
$K\subset 2^{[m]}$. We would assume that $K$ is not the $m$-simplex
$\D^m=2^{[m]}$, so $n<m$. If $I=\{i_1,\ldots,i_k\}\subset[m]$ is a simplex of
$K$, then we would write $I\in K$. By definition, the barycentric subdivision
of $K$ is the simplicial complex $\bs(K)$ whose simplices are chains
$I_1\subset I_2\subset\dots\subset I_p$ of embedded simplices of $K$.  In
particular, the vertices of $\bs(K)$ are in one to one correspondence with
simplices of $K$ of all dimensions. We denote the geometric realization of
$K$ (as a polyhedron) by $|K|$. Let $I^m$ be the standard unit
$m$-dimensional cube in $\R^m$:
$$
  I^m=\{(y_1,\ldots,y_m)\in\R^m\,:\;0\le y_i\le 1,\,i=1,\ldots,m\}.
$$
Every face of $I^m$ has the form
\begin{equation}
\label{ijface}
  F_{I\subset J}=\{(y_1,\ldots,y_m)\in I^m\: : \: y_i=0\text{ for }i\in
  I,\; y_j=1\text{ for }j\notin J\},
\end{equation}
were $I\subset J$ are two (possibly empty) subsets of $[m]$.  Now assign to
each subset $I=\{i_1,\ldots,i_k\}\subset[m]$ the vertex $v_I=F_{I\subset I}$
of the cube $I^m$. One has $v_I=(\e_1,\ldots,\e_m)$, where $\e_i=0$ if $i\in
I$ and $\e_i=1$ otherwise. Viewing $I$ as a vertex of the barycentric
subdivision $\bs(\D^m)$ of an $m$-simplex, we see that the correspondence
$I\mapsto v_I$ extends to a piecewise linear embedding $i_c$ of the
polyhedron $|\bs(\D^m)|$ into the (boundary complex of) unit cube $I^m$.
Under this embedding, the vertices of $\D^m$ are mapped to the vertices
$(1,\ldots,1,0,1,\ldots,1)$ of the cube $I^m$, while the point in the
interior of $|\D^m|$ (viewed as a vertex of $\bs(\D^m)$) is mapped to the
vertex $(0,\ldots,0)$ of $I^m$.  Hence, the whole image of $|\bs(\D^m)|$ is
the set of $m$ facets of $I^m$ meeting at the vertex $(0,\ldots,0)$.
Moreover, given a pair $I$, $J$ of non-empty subsets of $[m]$ such that
$I\subset J$, all simplices of $\bs(\D^m)$ of the form $I=I_1\subset
I_2\subset\dots\subset I_k=J$ are mapped to the same face $F_{I\subset J}$ of
$I^m$ (see~(\ref{ijface})).

Our simplicial complex $K^{n-1}$ can be viewed as a subcomplex in $\D^m$.
Hence, the above constructed map $i_c:|\bs(\D^m)|\to I^m$ embeds $|\bs(K)|$
piecewise linearly into the boundary complex of $I^m$. The image
$i_c(|\bs(K)|)$ is a certain $(n-1)$-dimensional cubical complex, which we
denote $\cub(K)$. Thus, we have proved the following statement.
\begin{proposition}
\label{cubk}
  There is a piecewise linear embedding $i_c$ of $|\bs(\D^m)|$ into the
  boundary complex of $I^m$ such that for any simplicial complex
  $K\subset\D^m$ on $m$ vertices the image $i_c(|\bs(K)|)=:\cub(K)$ is the
  union of faces~(\ref{ijface}) corresponding to all pairs $I\subset J$ of
  embedded simplices of $K$.\ep
\end{proposition}

\begin{remark}
Cubes of the cubical subdivision $\cub(K)$ of the polyhedron $|K|$ are formed
by simplices of the barycentric subdivision $\bs(K)$. This cubical
subdivision was employed in some previous researches for different purposes
(see, e.g.,~\cite[p.~434]{DJ}). The above cubical embedding $i_c:\cub(K)\to
I^m$ was used previously in~\cite{SS} to study which cubical complexes can be
embedded into the standard cubical lattice.
\end{remark}

The map $i_c:|\bs(\D^m)|\to I^m$ can be extended to a piecewise linear map
$\cone(i_c):|\cone(\bs(\D^m))|\to I^m$ by taking the vertex of the cone to
the vertex $(1,\ldots,1)$ of the cube $I^m$. (Note that the cone over the
barycentric subdivision of a $k$-simplex is identified with the standard
triangulation of a $(k+1)$-cube.) Now the image of $|\cone(\bs(\D^m))|$ is
the whole $I^m$, so $\cone(i_c)$ is a PL homeomorphism. The image $i_c\bigl(
|\cone(\bs(K))| \bigr)\subset I^m$ is another, this time $n$-dimensional,
cubical subcomplex of $I^m$, which we denote $\cc(K)$. This cubical complex
is explicitly described by the following proposition.
\begin{proposition}
\label{cck}
  For any simplicial complex $K$ on $m$ vertices there is a piecewise linear
  embedding of the polyhedron $|\cone(\bs(K))|$ into the boundary complex
  of $I^m$ such that its image $\cc(K)$ is the union of faces
  \begin{equation}
  \label{jface}
    F_J=\{(y_1,\ldots,y_m)\in I^m\: : \: y_j=1\text{ for }j\notin J\}
    \subset I^m
  \end{equation}
  and all their subfaces, where $J$ ranges over all simplices of $K$.\ep
\end{proposition}

According to~(\ref{ijface}), $F_J=F_{\varnothing\subset J}$. Hence, any
subface of $F_J$ is $F_{I\subset J}$ for some (possibly empty) $I\subset J$.
It follows that
$$
  \cub(K)=\mathop{\bigcup_{J\in K}}\limits_{I\ne\varnothing}F_{I\subset
  J},\quad \cc(K)=\bigcup_{J\in K}F_{I\subset J}.
$$
\begin{remark}
Viewed as a topological space, $|\cub(K)|$ is homeomorphic to $|K|$, while
$|\cc(K)|$ is homeomorphic to $|\cone(K)|$. Viewing $\cone(K)$ as a
simplicial complex, one may construct the cubical complex
$\cub\bigl(\cone(K)\bigr)$, which is also homeomorphic to $|\cone(K)|$.
However, as {\it cubical complexes}, $\cc(K)$ and $\cub\bigl(\cone(K)\bigr)$
differ.
\end{remark}

The cubical complex $\cc(K)$ was introduced in~\cite{BP2} and studied
in~\cite{BP3}, \cite{BP4} in connection with simple (and simplicial)
polytopes and subspace arrangements.
\begin{figure}
  \begin{picture}(120,45)
  \put(15,5){\line(1,0){25}}
  \put(15,5){\line(0,1){25}}
  \put(40,5){\line(1,1){10}}
  \put(50,15){\line(0,1){25}}
  \put(50,40){\line(-1,0){25}}
  \put(25,40){\line(-1,-1){10}}
  \put(40,5){\circle*{2}}
  \put(15,30){\circle*{2}}
  \put(50,40){\circle*{2}}
  \put(40,30){\line(1,1){10}}
  \multiput(25,15)(5,0){5}{\line(1,0){3}}
  \multiput(25,15)(0,5){5}{\line(0,1){3}}
  \multiput(25,15)(-5,-5){2}{\line(-1,-1){4}}
  \put(26,16){$0$}
  \put(22,-2){a)\ $K=\cdots$}
  \put(75,5){\line(1,0){25}}
  \put(75,5){\line(0,1){25}}
  \put(100,5){\line(1,1){10}}
  \put(110,15){\line(0,1){25}}
  \put(110,40){\line(-1,0){25}}
  \put(85,40){\line(-1,-1){10}}
  \put(100,5){\circle*{2}}
  \put(110,15){\circle*{2}}
  \put(75,30){\circle*{2}}
  \put(85,40){\circle*{2}}
  \put(75,5){\circle*{2}}
  \put(110,40){\circle*{2}}
  \put(100,30){\line(1,1){10}}
  \multiput(75,29.3)(0,0.1){16}{\line(1,1){10}}
  \multiput(100,4.3)(0,0.1){16}{\line(1,1){10}}
  \multiput(85,15)(5,0){5}{\line(1,0){3}}
  \multiput(85,15)(0,5){5}{\line(0,1){3}}
  \multiput(85,15)(-5,-5){2}{\line(-1,-1){4}}
  \put(86,16){$0$}
  \put(82,-2){b)\ $K=\bigtriangleup$}
  \put(15,30){\line(1,0){25}}
  \put(40,5){\line(0,1){25}}
  \put(75,30){\line(1,0){25}}
  \put(100,5){\line(0,1){25}}
  \linethickness{1mm}
  \put(75,5){\line(1,0){25}}
  \put(85,40){\line(1,0){25}}
  \put(75,5){\line(0,1){25}}
  \put(110,15){\line(0,1){25}}
  \end{picture}
  \caption{The cubical complex $\cub(K)$.}
  \end{figure}
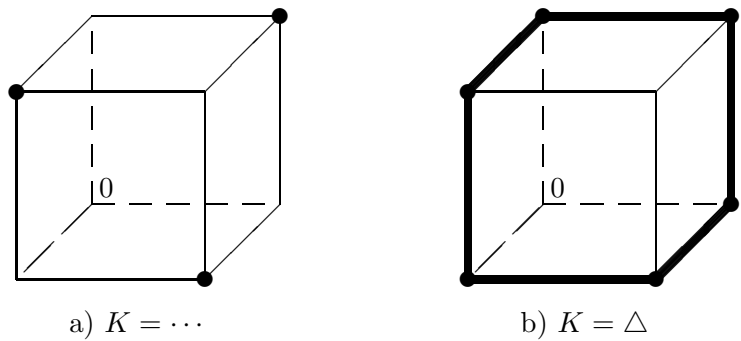
  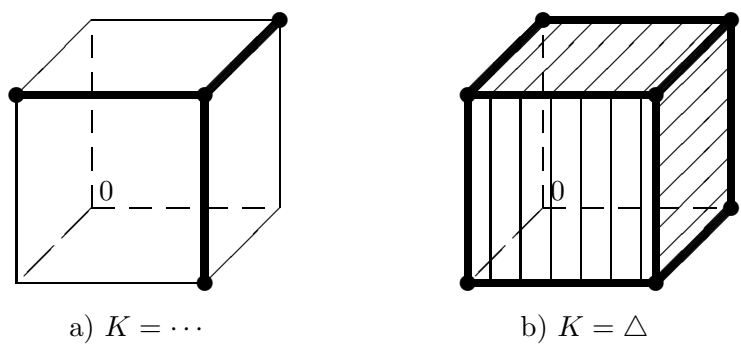
\begin{figure}
  \begin{picture}(120,45)
  \put(15,5){\line(1,0){25}}
  \put(15,5){\line(0,1){25}}
  \put(40,5){\line(1,1){10}}
  \put(50,15){\line(0,1){25}}
  \put(50,40){\line(-1,0){25}}
  \put(25,40){\line(-1,-1){10}}
  \put(40,5){\circle*{2}}
  \put(40,30){\circle*{2}}
  \put(15,30){\circle*{2}}
  \put(50,40){\circle*{2}}
  \multiput(40,29.3)(0,0.1){16}{\line(1,1){10}}
  \multiput(25,15)(5,0){5}{\line(1,0){3}}
  \multiput(25,15)(0,5){5}{\line(0,1){3}}
  \multiput(25,15)(-5,-5){2}{\line(-1,-1){4}}
  \put(26,16){$0$}
  \put(22,-2){a)\ $K=\cdots$}
  \put(75,5){\line(1,0){25}}
  \put(75,5){\line(0,1){25}}
  \put(100,5){\line(1,1){10}}
  \put(110,15){\line(0,1){25}}
  \put(110,40){\line(-1,0){25}}
  \put(85,40){\line(-1,-1){10}}
  \put(100,5){\circle*{2}}
  \put(110,15){\circle*{2}}
  \put(100,30){\circle*{2}}
  \put(75,30){\circle*{2}}
  \put(85,40){\circle*{2}}
  \put(75,5){\circle*{2}}
  \put(110,40){\circle*{2}}
  \multiput(100,29.3)(0,0.1){16}{\line(1,1){10}}
  \multiput(75,29.3)(0,0.1){16}{\line(1,1){10}}
  \multiput(100,4.3)(0,0.1){16}{\line(1,1){10}}
  \multiput(85,15)(5,0){5}{\line(1,0){3}}
  \multiput(85,15)(0,5){5}{\line(0,1){3}}
  \multiput(85,15)(-5,-5){2}{\line(-1,-1){4}}
  \multiput(78,5)(4,0){6}{\line(0,1){25}}
  \multiput(78,30)(4,0){6}{\line(1,1){10}}
  \multiput(100,8)(0,4){6}{\line(1,1){10}}
  \put(86,16){$0$}
  \put(82,-2){b)\ $K=\bigtriangleup$}
  \linethickness{1mm}
  \put(15,30){\line(1,0){25}}
  \put(40,5){\line(0,1){25}}
  \put(75,30){\line(1,0){25}}
  \put(75,5){\line(1,0){25}}
  \put(85,40){\line(1,0){25}}
  \put(100,5){\line(0,1){25}}
  \put(75,5){\line(0,1){25}}
  \put(110,15){\line(0,1){25}}
  \end{picture}
  \caption{The cubical complex $\cc(K)$.}
  \end{figure}
\begin{example}
  The cubical complex $\cub(K)$ for $K$ a disjoint union of
  3 vertices and the boundary complex of a 2-simplex is shown on
  Figure~1~a) and~b) respectively. The corresponding cubical complexes
  $\cc(K)$ are indicated on Figure~2~a) and~b).
\end{example}

\section{Equivariant moment-angle complexes}

Let $(D^2)^m$ denote the unit poly-disk in the complex space:
$$
  (D^2)^m=\{ (z_1,\ldots,z_m)\in\C^m:\: |z_i|\le1,\quad i=1,\ldots,m \}.
$$
The unit cube $I^m$ can be viewed as the orbit space for the standard action
of $m$-dimensional torus $T^m$ on $(D^2)^m$ by coordinatewise rotations. The
orbit map $\rho:(D^2)^m\to I^m$ can be given by $(z_1,\ldots,z_m)\to
(|z_1|^2,\ldots,|z_m|^2)$. For each face $F_{I\subset J}$ of $I^m$
(see~(\ref{ijface})) define
\begin{multline}
\label{ijblock}
  B_{I\subset J}:=\rho^{-1}(F_{I\subset J})\\ =\{(z_1,\ldots,z_m)\in
  (D^2)^m\: : \: z_i=0\text{ for }i\in I,\; |z_j|=1\text{ for }j\notin J\}.
\end{multline}
It follows that if $\#I=i$, $\#J=j$, then $B_{I\subset
J}\cong(D^2)^{j-i}\times T^{m-j}$, where disk factors $D^2\subset(D^2)^{j-i}$
correspond to elements from $J\setminus I$, while circle factors $S^1\subset
T^{m-j}$ correspond to elements from $[m]\setminus J$. Introducing the polar
coordinate system in $(D^2)^m$, we see that $B_{I\subset J}$ is parametrized
by $(j-i)$ radial (or moment) and $(m-i)$ angle coordinates. Here we come to
the following definition.

\begin{definition}
\label{ma}
  Let $C$ be a cubical subcomplex of $I^m$.  The {\it moment-angle complex}
  $\ma(C)$ corresponding to $C$ is the $T^m$-invariant decomposition of
  $\rho^{-1}(|C|)$ to ``moment-angle" blocks $B_{I\subset J}$
  (see~(\ref{ijblock})) corresponding to faces $F_{I\subset
  J}\subset|C|\subset I^m$. Hence, $\ma(C)$ is defined from the commutative
  diagram
  $$
  \begin{CD}
    \ma(C) @>>> (D^2)^m\\
    @VVV @VV\rho V\\
    |C| @>>> I^m
  \end{CD}.
  $$
\end{definition}
It follows that the torus $T^m$ acts on $\ma(C)$ with orbit space $|C|$.

The moment-angle complexes corresponding to the introduced above cubical
complexes $\cub(K)$ and $\cc(K)$ (see propositions~\ref{cubk} and~\ref{cck})
will be denoted $\wk$ and $\zk$ respectively. Thus, we have
\begin{align}
\label{zkwk}
  \begin{CD}
    \wk @>>> (D^2)^m\\
    @VVV @VV\rho V\\
    |\cub(K)| @>>> I^m
  \end{CD} &&
  \text{and} &&
  \begin{CD}
    \zk @>>> (D^2)^m\\
    @VVV @VV\rho V\\
    |\cc(K)| @>>> I^m
  \end{CD}, &
\end{align}
where horizontal arrows are embeddings, while vertical ones are orbit maps
for certain $T^m$-actions. It follows that $\dim\zk=m+n$ and $\dim\wk=m+n-1$.

\begin{figure}
  \begin{picture}(120,40)
  \put(60,15){\oval(30,30)}
  \put(60,15){\line(1,0){15}}
  \put(60,15){\circle*{1.5}}
  \put(75,15){\circle*{1.5}}
  \put(61,16){0}
  \put(76,16){1}
  \put(68,16){$I$}
  \put(42,16){$T$}
  \put(52,20){$D$}
  \end{picture}
  \caption{Cellular decomposition of $D^2$.}
\end{figure}
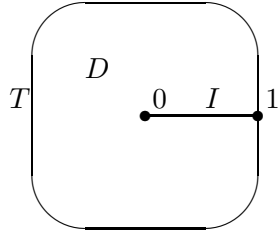
Let us consider the cellular decomposition of $D^2$ with one 2-cell $D$, two
1-cells $I$, $T$, and two 0-cells 0, 1 (see Figure~3). It defines a
$T^m$-invariant cellular decomposition of the poly-disk $(D^2)^m$ with $5^m$
cells. Each cell of this decomposition is the product of cells of 5 different
types: $D_i$, $I_i$, $0_i$, $T_i$, and $1_i$, $i=1,\ldots,m$. We will encode
cells of $(D^2)^m$ by ``words" of type $D_II_J0_LT_P1_Q$, where
$I,J,L,P,Q$ are pairwise disjoint subsets of $[m]$ such that $I\cup J\cup
L\cup P\cup Q=[m]$. Sometimes we would drop the last factor $1_Q$, so in our
notations $D_II_J0_LT_P=D_II_J0_LT_P1_{[m]\setminus I\cup J\cup L\cup P}$. It
follows that the closure of $D_II_J0_LT_P1_Q$ is homeomorphic to the product
of $\#I$ disks, $\#J$ segments, and $\#P$ circles. The constructed cellular
decomposition of $(D^2)^m$ allows to identify moment-angle complexes as
certain cellular subcomplexes in $(D^2)^m$.
\begin{lemma}
\label{macell}
  For any cubical subcomplex $C$ of $I^m$ the corresponding
  moment-angle complex $\ma(C)$ is a $T^m$-invariant cellular
  subcomplex of $(D^2)^m$.
\end{lemma}
\begin{proof}
Since $\ma(C)$ is a union of ``moment-angle" blocks
$B_{I\subset J}$ (see~(\ref{ijblock})), and each $B_{I\subset J}$ is
obviously $T^m$-invariant, the whole $\ma(C)$ is also
$T^m$-invariant. In order to show that $\ma(C)$ is a cellular
subcomplex of $(D^2)^m$ (with respect to the above constructed cellular
decomposition) we just mention that $B_{I\subset J}$ is the closure of
cell $D_{J\setminus I}I_\varnothing0_IT_{[m]\setminus J}1_\varnothing$.
\end{proof}

\section{Cohomology of $\zk$, subspace arrangements, and numbers of faces
of~$K$}

Here we study the moment-angle complex $\zk$ corresponding to the cubical
complex $\cc(K)\subset I^m$ (see~(\ref{zkwk})). Remember that $K$ is an
$(n-1)$-dimensional simplicial complex on $m$ vertices, and $\cc(K)$
topologically is the cone over $K$. By definition, $\zk$ is
a union of certain moment-angle blocks $B_{I\subset J}\subset(D^2)^m$ with
$I=\varnothing$ (see Proposition~\ref{cck}). In analogy with~(\ref{jface}),
we put
$$
  B_J:=B_{\varnothing\subset J}= \rho^{-1}(F_J)=\{(z_1,\ldots,z_m)\in
  (D^2)^m\: : \: |z_j|=1\text{ for }j\notin J\}.
$$
Hence, $\zk=\bigcup_{J\in K}B_J$, where $B_J\cong(D^2)^j\times T^{m-j}$,
$j=\#J$. This remark allows to simplify the cellular decomposition
constructed in the previous section for $\ma(C)$ in the case when
$\ma(C)=\zk$. To do this we replace the union of cells $0$, $I$, $D$ (see
Figure~3) by one 2-dimensional cell. To simplify notations we
denote this 2-cell by $D$ throughout this section. The resulting
$T^m$-invariant cellular decomposition of $(D^2)^m$ now has $3^m$ cells, each
of which is the product of 3 different types of cells: $D_i$, $T_i$, and
$1_i$, $i=1,\ldots,m$. In this section we encode these cells of $(D^2)^m$ as
$D_IT_P1_Q$, where $I,P,Q$ are pairwise disjoint subsets of $[m]$ such that
$I\cup P\cup Q=[m]$. Usually we would denote the cell $D_IT_P1_Q$ just by
$D_IT_P$, so $D_IT_P=D_IT_P1_{[m]\setminus I\cup P}$.  Since
$B_J=B_{\varnothing\subset J}$ is the closure of cell $D_JT_{[m]\setminus
J}1_\varnothing$, the moment-angle complex $\zk$ is a $T^m$-invariant
cellular subcomplex of $(D^2)^m$ (with respect to the new $3^m$-cell
decomposition). The complex $\zk$ consists of all cells
$D_IT_P\subset(D^2)^m$ such that $I$ is a simplex of~$K$.

\begin{remark}
Note that for general $C$ the moment-angle complex $\ma(C)$ is {\it not} a
cellular subcomplex for the $3^m$-cell decomposition of $(D^2)^m$.
\end{remark}

The cohomology ring of $\zk$ was described in~\cite{BP2},~\cite{BP3} (in the
case when $K$ is a polytopal sphere) and in~\cite{BP4} (for general $K$).
Before going further, we review some results of these papers.

Throughout the rest of this paper we work over some field $\k$, referred to
as the ground field. Let $\k[v_1,\ldots,v_m]$ be the polynomial algebra, and
$\Lambda[u_1,\ldots,u_m]$ the exterior algebra over $\k$ on $m$ generators.
We make both algebras graded by putting $\deg(v_i)=2$, $\deg(u_i)=1$.

\begin{definition}
\label{fr}
  The {\it face ring\/} (or the {\it Stanley--Reisner ring\/}) $\k(K)$ of
  simplicial complex $K$ is the quotient ring $\k[v_1,\ldots,v_m]/I$, where
  the ideal $I$ is generated by all square-free monomials $v_{i_1}\cdots
  v_{i_s}$, $1\le i_1<\dots<i_s\le m$, such that $\{i_1,\ldots,i_s\}$ is not
  a simplex of K.
\end{definition}

For any subset $I=\{i_1,\ldots,i_k\}\in[m]$ denote by $L_I$ the coordinate
plane in $\C^m$ consisting of points whose $i_1,\ldots,i_k$ coordinates
vanish:
\begin{equation}
\label{li}
  L_I=\{(z_1,\ldots,z_m)\in\C^m\::\:z_{i_1}=\cdots=z_{i_k}=0\}.
\end{equation}
Each simplicial complex $K$ on $m$ vertices defines a {\it complex coordinate
subspace arrangement} $\mathcal A(K)$. The latter is the set of all planes
$L_I$ such that $I$ is not a simplex of $K$:
$$
  \mathcal A(K)=\{L_I\::\:I\notin K\}.
$$
Define the {\it support} of $\mathcal A(K)$ as $|\mathcal
A(K)|=\bigcup_{I\notin K}L_I\subset\C^m$ and the {\it complement}
$U(K)=\C^m\setminus|\mathcal A(K)|$, that is
\begin{equation}
\label{compl}
  U(K)=\C^m\setminus\bigcup_{I\notin K}L_I.
\end{equation}
It can be easily seen that the complement of {\it any} coordinate subspace
arrangement in $\C^m$ (i.e., the complement of any set of planes~(\ref{li}))
is $U(K)$ for some $K$. Note that $U(K)$ is invariant with respect to the
standard $T^m$-action on $\C^m$. The following lemma establishes the
connection between moment-angle complexes and complements of coordinate
subspace arrangements.
\begin{lemma}
\label{he}
  The complement $U(K)$ is $T^m$-equivariantly homotopy equivalent to the
  moment-angle complex $\zk$.
\end{lemma}
\begin{proof}
See~\cite[Lemma~2.13]{BP4}
\end{proof}
The next theorem describes the cohomology algebra of $\zk$.
\begin{theorem}
\label{arrang}
  The following isomorphisms of graded algebras holds:
  \begin{equation}
  \label{toralg}
    H^{*}\bigl(\zk\bigr)\cong\Tor_{\k[v_1,\ldots,v_m]}\bigl(\k(K),\k\bigr)
    \cong H\bigl[\k(K)\otimes\Lambda[u_1,\ldots,u_m],d\bigr],
  \end{equation}
  where $\k(K)=\k[v_1,\ldots,v_m]/I$ is the face ring, and the differential
  $d$ is defined by $d(v_i)=0$, $d(u_i)=v_i$, $i=1,\ldots,m$.
\end{theorem}
\begin{proof}
See~\cite[theorems 3.2 and 3.3]{BP4}
\end{proof}
Due to Lemma~\ref{he}, isomorphism (\ref{toralg}) also holds for the
cohomology algebra of the complement of a coordinate subspace arrangement in
$\C^m$. The first isomorphism of~(\ref{toralg}) is proved by applying the
Eilenberg--Moore spectral sequence to some $T^m$-bundles. The second
isomorphism follows from the Koszul complex for the
$\k[v_1,\ldots,v_m]$-module $\k(K)$.

\begin{remark}
As we mentioned in the introduction, the Betti numbers of the complement of a
{\it real\/} coordinate subspace arrangement were calculated in terms of
resolution of the Stanley--Reisner ring in~\cite{GPW}. The latter paper also
contains the formulation of the multiplicative isomorphism~(\ref{toralg}) for
complex coordinate subspace arrangements (see~\cite[Thm.~3.6]{GPW}) with a
reference to yet unpublished paper by Babson and Chan. We note also, that, as
it was observed in~\cite{GPW}, there is no isomorphism between the cohomology
algebra of a real coordinate subspace arrangement and the corresponding
Stanley--Reisner $\Tor$-algebra.
\end{remark}

The $\Tor$-algebra from~(\ref{toralg}) is naturally a bigraded algebra with
$\bideg(v_i)=(0,2)$, $\bideg(u_i)=(-1,2)$, and differential $d$ adding
$(1,0)$ to bidegree. Since the differential does not change the second
grading, the whole algebra is decomposed into the sum of differential
algebras consisting of elements with fixed second degree. Below we deduce
some important consequences of this bigraded structure.

\begin{remark}
Note that according to our agreement the first degree in the $\Tor$-algebra
is {\it non-positive}. This corresponds to numerating the terms of Koszul
$\k[v_1,\ldots,v_m]$-free resolution of $\k$ by non-positive integers. In
such notations the Koszul complex
$[\k(K)\otimes\Lambda[u_1,\ldots,u_m],d\bigr]$ becomes a {\it cochain}
complex, and $\Tor_{\k[v_1,\ldots,v_m]}\bigl(\k(K),\k\bigr)$ is its {\it
cohomology}, not the homology as usually regarded. This is the standard trick
used for applying the Eilenberg--Moore spectral sequence, see~\cite{Sm}. It
also explains why we write
$\Tor^{*,*}_{\k[v_1,\ldots,v_m]}\bigl(\k(K),\k\bigr)$ instead of usual
$\Tor_{*,*}^{\k[v_1,\ldots,v_m]}\bigl(\k(K),\k\bigr)$.
\end{remark}

Following~\cite{BP2}, define the subcomplex $\mathcal C^{*}(K)$ of the
cochain complex $\bigl[\k(K)\otimes\Lambda[u_1,\ldots,u_m],d\bigr]$
(see~(\ref{toralg})) as follows. As a $\k$-module, $\mathcal
C^{*}(K)=\bigoplus_{q=0}^{m}\mathcal C^{-q}(K)$, where $\mathcal C^{-q}(K)$
is generated by monomials $u_{j_1}\ldots u_{j_q}$ and $v_{i_1}\ldots
v_{i_p}u_{j_1}\ldots u_{j_q}$ such that $\{i_1,\ldots,i_p\}$ is a simplex of
$K$ and $\{i_1,\ldots,i_p\}\cap\{j_1,\ldots,j_q\}=\varnothing$. Since
$d(u_i)=v_i$, we have $d\bigl(\mathcal C^{-q}(K)\bigr)\subset\mathcal
C^{-q+1}(K)$ and, therefore, $\mathcal C^{*}(K)$ is a cochain subcomplex.
Moreover, $\mathcal C^{*}(K)$ inherits the bigraded module structure from
$\k(K)\otimes\Lambda[u_1,\ldots,u_m]$, with differential $d$ adding $(1,0)$
to bidegree. Hence, we have the additive inclusion (i.e., the monomorphism of
bigraded modules) $i:\mathcal
C^{*}(K)\hookrightarrow\k(K)\otimes\Lambda[u_1,\ldots,u_m]$. Finally,
$\mathcal C^{*}(K)$ can be viewed as an algebra in obvious way, however this
time this is not a subalgebra of $\k(K)\otimes\Lambda[u_1,\ldots,u_m]$
(since, for instance, $v_1^2=0$ in $\mathcal C^{*}(K)$ but not in
$\k(K)\otimes\Lambda[u_1,\ldots,u_m]$). However, we have the multiplicative
{\it projection} (i.e., the epimorphism of bigraded algebras)
$j:\k(K)\otimes\Lambda[u_1,\ldots,u_m]\to\mathcal C^{*}(K)$. The additive
inclusion $i$ and the multiplicative projection $j$ obviously satisfy $j\cdot
i=id$.

\begin{lemma}
\label{iscoh}
  Cochain complexes $\bigl[\k(K)\otimes\Lambda[u_1,\ldots,u_m],d\bigr]$
  and $[\mathcal C^{*}(K),d]$ have same cohomology. Hence, the following
  isomorphism of bigraded $\k$-modules holds:
  $$
    H[\mathcal
    C^{*}(K),d]\cong\Tor^{*}_{k[v_1,\ldots,v_m]}\bigl(\k(K),\k\bigr).
  $$
\end{lemma}
\begin{proof}See~\cite[Lemma~5.3]{BP2}.
\end{proof}

In the sequel we denote (square-free) monomials $v_{i_1}\ldots
v_{i_p}u_{j_1}\ldots u_{j_q}\in\k(K)\otimes\Lambda[u_1,\ldots,u_m]$ by
$v_Iu_J$, where $I=\{i_1,\ldots,i_p\}$, $J=\{j_1,\ldots,j_q\}$ are
multiindices.

Now we recall our cellular decomposition of $\zk$, whose cells are $D_IT_J$,
where $I,J\subset[m]$, $I$ is a simplex of $K$, and $I\cap
J=\varnothing$. Let $\mathcal C_{*}(\zk)$ and $\mathcal C^{*}(\zk)$ denote
the corresponding cellular chain and cochain complex respectively. Both
complexes (or differential algebras) $\mathcal C^{*}(\zk)$ and $\mathcal
C^{*}(K)$ have same cohomology $H^{*}(\zk)$. Cochain complex $\mathcal
C^{*}(\zk)$ has basis consisting of elements $(D_IT_J)^{*}$ dual to
$D_IT_J\in\mathcal C_{*}(\zk)$ (the latter is viewed as a cellular chain).
Note that the cochain algebra $\mathcal C^{*}(\zk)$ is multiplicatively
generated by the elements $T_i^{*}$, $D_i^{*}$, $i,j=1,\ldots,m$, (of
dimension 1 and 2 respectively), while $\mathcal C^{*}(K)$ is
multiplicatively generated by $u_i$, $v_i$, $i,j=1,\ldots,m$. The following
theorem shows that these two algebras are the same.

\begin{theorem} \label{cellcom}
  The correspondence $v_Iu_J\mapsto(D_IT_J)^{*}$ establishes a canonical
  isomorphism of differential graded algebras $\mathcal C^{*}(K)$ and
  $\mathcal C^{*}(\zk)$.
\end{theorem}
\begin{proof}
It follows directly from the definitions of $\mathcal C^{*}(K)$ and $\mathcal
C^{*}(\zk)$ that the constructed map is an isomorphism of graded algebras.
So, it remains to prove that it commutes with differentials. Let $d$, $d_c$
and $\partial_c$ denote the differentials in $\mathcal C^{*}(K)$, $\mathcal
C^{*}(\zk)$ and $\mathcal C_{*}(\zk)$ respectively.  Since $d(v_i)=0$,
$d(u_i)=v_i$, we need to show that $d_c(D_i^{*})=0$, $d_c(T_i^{*})=D_i^{*}$.
We have $\partial_c(D_i)=T_i$, $\partial_c(T_i)=0$. Since any 2-cell of $\zk$
is either $D_j$ or $T_{jk}$, $k\ne j$, it follows that
$$
  (d_cT_i^{*},D_j)=(T_i^{*},\partial_cD_j)=(T_i^{*},T_j)=\delta_{ij},\quad
  (d_cT_i^{*},T_{jk})=(T_i^{*},\partial_cT_{jk})=0,
$$
where $\delta_{ij}=1$ if $i=j$ and $\delta_{ij}=0$ otherwise. Hence,
$d(T_i^{*})=D_i^{*}$. Further, since any 3-cell of $\zk$
is either $D_jT_k$ or $T_{j_1j_2j_3}$, it follows that
\begin{align*}
  &(d_cD_i^{*},D_jT_k)=(D_i^{*},\partial_c(D_jT_k))=
  (D_i^{*},T_{jk})=0,\\
  &(d_cD_i^{*},T_{j_1j_2j_3})=(D_i^{*},\partial_cT_{j_1j_2j_3})=0.
\end{align*}
Hence, $d_c(D_i^{*})=0$.
\end{proof}
In the sequel we would not distinguish cochain complexes $\mathcal C^{*}(K)$
and $\mathcal C^{*}(\zk)$ and their elements $u_i$ and $T_i^{*}$, $v_i$ and
$D_i^{*}$. The above theorem provides two methods for calculating the
(co)homology of $\zk$: either by means of the differential (bi)graded algebra
$[\mathcal C^{*}(K),d]$, where $\mathcal
C^{*}(K)\subset\k(K)\otimes\Lambda[u_1,\ldots,u_m]$ (as modules), or using
the cellular chain complex $[\mathcal C_{*}(\zk),\partial_c]$.

Now we recall that the algebra $[\mathcal C^{*}(K),d]$ is bigraded.
Theorem~\ref{cellcom} shows that the cellular chain complex $[\mathcal
C_{*}(\zk),\partial_c]$ can be also made bigraded by setting
\begin{equation}\label{bgcellz}
  \bideg(D_i)=(0,2),\quad \bideg(T_i)=(-1,2),\quad \bideg(1_i)=(0,0).
\end{equation}
The differential $\partial_c$ adds $(-1,0)$ to bidegree, and the cellular
homology of $\zk$ also acquires a bigraded structure. Let us assume now that
the ground field $\k$ is of zero characteristic (e.g., $\k=\Q$ is the field
of rational numbers). Define the bigraded Betti numbers
\begin{equation}
\label{bbn}
  b_{-q,2p}(\zk)=\dim H_{-q,2p}[\mathcal C_{*}(\zk),\partial_c],
  \quad q,p=0,\ldots,m.
\end{equation}
Theorem~\ref{cellcom} and Lemma~\ref{iscoh} show that
$$
  b_{-q,2p}(\zk)=\dim\Tor^{-q,2p}_{\k[v_1,\ldots,v_m]}\bigl(\k(K),\k\bigr)
$$
(i.e., $b_{-q,2p}(\zk)$ is the dimension of $(-q,2p)$-th bigraded component
of the cohomology algebra
$H\bigl[\k(K)\otimes\Lambda[u_1,\ldots,u_m],d\bigr]$). For the ordinary Betti
numbers $b_k(\zk)$ holds
$$
  b_k(\zk)=\sum_{-q+2p=k}b_{-q,2p}(\zk),\quad k=0,\ldots,m+n.
$$

Below we describe some basic properties of bigraded Betti
numbers~(\ref{bbn}).
\begin{lemma}\label{bbgen}
  Let $K^{n-1}$ be an $(n-1)$-dimensional simplicial complex with $m=f_0$
  vertices and $f_1$ edges, and let $\zk$ be the corresponding moment-angle
  complex, $\dim\zk=m+n$. Then
  \begin{itemize}
  \item[\it (a)] $b_{0,0}(\zk)=b_0(\zk)=1$,\quad $b_{0,2p}(\zk)=0$ if $p>0$;
  \item[\it (b)] $b_{-q,2p}=0$ if $p>m$ or $q>p$;
  \item[\it (c)] $b_1(\zk)=b_2(\zk)=0$;
  \item[\it (d)] $b_3(\zk)=b_{-1,4}(\zk)=\binom{f_0}2-f_1$;
  \item[\it (e)] $b_{-q,2p}(\zk)=0$ if $q\ge p>0$ or $p-q>n$;
  \item[\it (f)] $b_{m+n}(\zk)=b_{-(m-n),2m}(\zk)$.
  \end{itemize}
\end{lemma}
\begin{proof}
We use the cochain complex $\mathcal
C^{*}(K)\subset\k(K)\otimes\Lambda[u_1,\ldots,u_m]$ for calculations. The
module $\mathcal C^{*}(K)$ has basis consisting of monomials $v_Iu_J$ with
$I\in K$ and $I\cap J=\varnothing$. Since $\bideg v_i=(0,2)$, $\bideg
u_j=(-1,2)$, the bigraded component $\mathcal C^{-q,2p}(K)$ is generated by
monomials $v_Iu_J$ with $\#I=p-q$ and $\#J=q$. In particular, $\mathcal
C^{-q,2p}(K)=0$ if $p>m$ or $q>p$, whence the assertion~(b) follows. To
prove~(a) we mention that $\mathcal C^{0,0}(K)$ is generated by 1, while any
$v_I\in\mathcal C^{0,2p}(K)$, $p>0$, is a coboundary, hence,
$H^{-q,2p}(K)=0$, $p>0$.

Now we are going to prove the assertion (e).  Since any $v_Iu_J\in\mathcal
C^{-q,2p}(K)$ has $I\in K$, and any simplex of $K$ is at most
$(n-1)$-dimensional, it follows that $\mathcal C^{-q,2p}(K)=0$ for $p-q>n$.
It follows from~(b) that $b_{-q,2p}(\zk)=0$ for $q>p$, so it remains to prove
that $b_{-q,2q}(\zk)=0$ for $q>0$. The module $\mathcal C^{-q,2q}(K)$ is
generated by monomials $u_J$, $\#J=q$. Since $d(u_i)=v_i$, it follows easily
that there no cocycles in $\mathcal C^{-q,2q}(K)$. Hence, $H^{-q,2q}(\zk)=0$.

To prove (c) we mention that $H^1(\zk)=H^{-1,2}(K)$ and
$H^{2}(\zk)=H^{-2,4}(\zk)$ (this follows from (a) and (b)). Hence, the
assertion~(c) follows from~(e).

As for (d), it follows from (a), (b) and (e) that $H^{3}(\zk)=H^{-1,4}(\zk)$.
The module $\mathcal C^{-1,4}(K)$ is generated by monomials $v_iu_j$, $i\ne
j$.  We have $d(v_iu_j)=v_iv_j$ and $d(u_iu_j)=v_iu_j-v_ju_i$. It follows
that $v_iu_j$ is a cocycle if and only if $\{i,j\}$ is not a 1-simplex in $K$;
in this case two cocycles $v_iu_j$ and $v_ju_i$ are cohomological. The
assertion~(d) now follows easily.

The remaining assertion (f) follows from the fact that the monomial
$u_Iv_J\in\mathcal C^{*}(K)$ of maximal total degree $m+n$ necessarily has
$\#I+\#J=m$, $\#J=n$, $\#I=m-n$.
\end{proof}
The above lemma shows that non-zero bigraded Betti numbers $b_{r,2p}(\zk)$,
$r\ne0$ appear only in the ``strip" bounded by the lines $r=-(m-1)$, $r=-1$,
$p+r=1$ and $p+r=n$ in the second quadrant (see Figure~4~(a)).
\begin{figure}
  \begin{picture}(115,60)
  \multiput(45,10)(0,5){10}{\line(-1,0){47}}
  \multiput(45,10)(-5,0){10}{\line(0,1){47}}
  \put(46,11){\small 0}
  \put(46,16){\small 2}
  \put(46,21){\small 4}
  \put(46,31){\vdots}
  \put(46,51){\small $2m$}
  \put(42,7){\footnotesize 0}
  \put(35.5,7){\footnotesize $-1$}
  \put(30,7){\small $\cdots$}
  \put(22.5,10){\line(0,-1){2}}
  \put(15,6){\footnotesize $-(m-n)$}
  \put(-1,7){\footnotesize $-m$}
  \put(41,11){\Large $*$}
  \multiput(36,21)(-5,5){7}{\Large $*$}
  \multiput(36,26)(-5,5){6}{\Large $*$}
  \multiput(36,31)(-5,5){5}{\Large $*$}
  \multiput(36,36)(-5,5){4}{\Large $*$}
  \put(7,0){a)\ arbitrary $K^{n-1}$}
  \multiput(110,10)(0,5){10}{\line(-1,0){47}}
  \multiput(110,10)(-5,0){10}{\line(0,1){47}}
  \put(111,11){\small 0}
  \put(111,16){\small 2}
  \put(111,21){\small 4}
  \put(111,31){\vdots}
  \put(111,51){\small $2m$}
  \put(107,7){\footnotesize 0}
  \put(100.5,7){\footnotesize $-1$}
  \put(95,7){\small $\cdots$}
  \put(87.5,10){\line(0,-1){2}}
  \put(80,6){\footnotesize $-(m-n)$}
  \put(64,7){\footnotesize $-m$}
  \put(106,11){\Large $*$}
  \multiput(101,21)(-5,5){3}{\Large $*$}
  \multiput(101,26)(-5,5){3}{\Large $*$}
  \multiput(101,31)(-5,5){3}{\Large $*$}
  \put(86,51){\Large $*$}
  \put(77,0){b)\ $|K|=S^{n-1}$}
  \end{picture}
  \caption{Possible locations of non-zero bigraded Betti numbers
  $b_{-q,2p}(\zk)$ are marked by $*$.}
\end{figure}
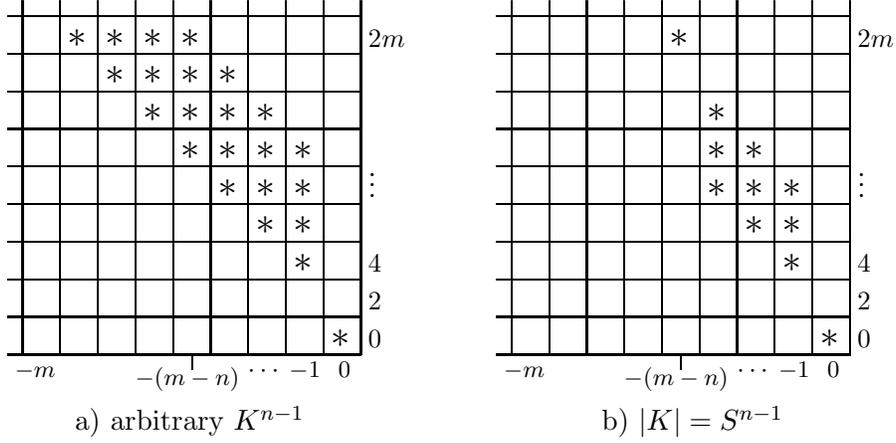

Let us consider now the bigraded cellular chain complex $[\mathcal
C_{*,*}(\zk),\partial_c]$. The homogeneous component $\mathcal C_{-q,2p}(\zk)$
consists of cellular chains $D_IT_J$ with $I\in K$, $\#I=p-q$, $\#J=q$. It
follows that
\begin{equation}\label{dimc}
  \dim\mathcal C_{-q,2p}(\zk)=f_{p-q-1}\binom{m-p+q}q
\end{equation}
(with usual agreement $\binom ij=0$ if $i<j$ or $j<0$), where $f_i$ is the
number of $i$-simplices of $K^{n-1}$ and $f_{-1}=1$. Since the differential
$\partial_c$ does not change the second degree, i.e.
$$
  \partial_c:\mathcal C_{-q,2p}(\zk)\to\mathcal C_{-q-1,2p}(\zk),
$$
the chain complex $\mathcal C_{*,*}(\zk)$ splits into the sum of chain
complexes as follows:
$$
  [\mathcal C_{*,*}(\zk),\partial_c]=\bigoplus_{p=0}^m
  [\mathcal C_{*,2p}(\zk),\partial_c].
$$
The similar decomposition holds also for the cellular cochain complex
$[\mathcal C^{*,*}(\zk),d_c]\cong[\mathcal C^{*,*}(K),d]$. Let $\chi_p(\zk)$
denote the Euler characteristic of complex
$[\mathcal C_{*,2p}(\zk),\partial_c]$, i.e.
\begin{equation}\label{chip}
  \chi_p(\zk)=\sum_{q=0}^m(-1)^q\dim\mathcal C_{-q,2p}(\zk)
  =\sum_{q=0}^m(-1)^qb_{-q,2p}(\zk).
\end{equation}
Define the generating polynomial $\chi(\zk,t)$ as
$$
  \chi(\zk,t)=\sum_{p=0}^m\chi_p(\zk)t^{2p}.
$$
The following theorem calculates this polynomial in terms of the numbers of
faces of $K$. It was firstly proved in~\cite{BP2} in the particular case of
polytopal~$K$.

\begin{theorem}
\label{gpz}
  For any $(n-1)$-dimensional simplicial complex $K$ with $m$ vertices holds
  \begin{equation}
  \label{hchi}
    \chi(\zk,t)=(1-t^2)^{m-n}(h_0+h_1t^2+\cdots+h_nt^{2n}),
  \end{equation}
  where $(h_0,h_1,\ldots,h_n)$ is the $h$-vector of $K$
  (see~{\rm (\ref{hvector})}).
\end{theorem}
\begin{proof}
It follows from (\ref{chip}) and (\ref{dimc}) that
\begin{equation}
\label{chipzk}
  \chi_p(\zk)=\sum_{j=0}^m(-1)^{p-j}f_{j-1}\binom{m-j}{p-j},
\end{equation}
Then
\begin{multline}
\label{chidir}
  \chi(\zk,t)=\sum_{p=0}^m\chi_p(K)t^{2p}=
  \sum_{p=0}^m\sum_{j=0}^mt^{2j}t^{2(p-j)}(-1)^{p-j}f_{j-1}
  \binom{m-j}{p-j}\\
  =\sum_{j=0}^mf_{j-1}t^{2j}(1-t^2)^{m-j}=
  (1-t^2)^m\sum_{j=0}^nf_{j-1}(t^{-2}-1)^{-j}.
\end{multline}
Denote $h(t)=h_0+h_1t+\cdots+h_nt^n$. Then it follows from~(\ref{hvector})
that
$$
  t^nh(t^{-1})=(t-1)^n\sum_{i=0}^nf_{i-1}(t-1)^{-i}.
$$
Substituting above $t^{-2}$ for $t$, we finally obtain from~(\ref{chidir})
$$
  \frac{\chi(\zk,t)}{(1-t^2)^m}=\frac{t^{-2n}h(t^2)}{(t^{-2}-1)^n}=
  \frac{h(t^2)}{(1-t^2)^n},
$$
which is equivalent to (\ref{hchi}).
\end{proof}
The above theorem allows to express the numbers of faces of a simplicial
complex in terms of bigraded Betti numbers of the corresponding moment-angle
complex $\zk$. The first important corollary of this is as follows.
\begin{corollary}\label{chizk}
  For any simplicial complex $K$ the Euler number of the corresponding
  moment-angle complex $\zk$ is zero.
\end{corollary}
\begin{proof}
We have
$$
  \chi(\zk)=\sum_{p,q=0}^m(-1)^{-q+2p}b_{-q,2p}(\zk)=
  \sum_{p=0}^m\chi_p(\zk)=\chi(\zk,1)
$$
Now the statement follows from (\ref{hchi}).
\end{proof}

\begin{remark}
Another way to prove the above corollary is to mention that the diagonal
subgroup $S^1\subset T^m$ always acts freely on the moment-angle complex
$\zk$ (see~\cite{BP2}). Hence, there exists a principal $S^1$-bundle
$\zk\to\zk/S^1$, which implies $\chi(\zk)=0$.
\end{remark}

\begin{corollary}\label{chiar}
  The Euler number of the complement of a complex coordinate subspace
  arrangement is zero.
\end{corollary}
\begin{proof}
This follows from the previous corollary and Lemma~\ref{he}.
\end{proof}

By definition (see Proposition~\ref{cck}), the cubical complex $\cc(K)$
always contains the vertex $(1,\ldots,1)\in I^m$. Hence, the torus
$T^m=\rho^{-1}(1,\ldots,1)$ is contained in $\zk$. Here $\rho:(D^2)^m\to I^m$
is the orbit map for the $T^m$-action (see~(\ref{zkwk})).

\begin{lemma}\label{torus}
  The inclusion $T^m=\rho^{-1}(1,\ldots,1)\hookrightarrow\zk$ is a cellular
  map homotopical to the map to a point, i.e.  the torus
  $T^m=\rho^{-1}(1,\ldots,1)$ is a contractible cellular subcomplex of $\zk$.
\end{lemma}
\begin{proof}
To prove that $T^m=\rho^{-1}(1,\ldots,1)$ is a cellular subcomplex of $\zk$
we just mention that this $T^m$ is the closure of the $m$-cell $D_\varnothing
T_{[m]}\subset\zk$. So, it remains to prove that $T^m$ is contractible within
$\zk$. To do this we show that the embedding $T^m\subset(D^2)^m$ is
homotopical to the map to the point $(1,\ldots,1)\in T^m\subset(D^2)^m$.  On
the first step we note that $\zk$ contains the cell $D_1T_{2,\ldots,m}$,
whose closure contains $T^m$ and is homeomorphic to $D^2\times T^{m-1}$.
Hence, our $T^m$ can be contracted to $1\times T^{m-1}$ within $\zk$. On the
second step we note that $\zk$ contains the cell $D_2T_{3,\ldots,m}$, whose
closure contains $1\times T^{m-1}$ and is homeomorphic to $D^2\times
T^{m-2}$. Hence, $1\times T^{m-1}$ can be contracted to $1\times 1\times
T^{m-2}$ within $\zk$, and so on. On the $k$th step we note that $\zk$
contains the cell $D_kT_{k+1,\ldots,m}$, whose closure contains
$1\times\cdots\times1\times T^{m-k+1}$ and is homeomorphic to $D^2\times
T^{m-k}$. Hence, $1\times\cdots\times1\times T^{m-k+1}$ can be contracted to
$1\times\cdots\times1\times T^{m-k}$ within $\zk$. We end up at the point
$1\times\cdots\times1$ to which the whole torus $T^m$ can be contracted.
\end{proof}

\begin{corollary}
  For any simplicial complex $K$ the moment-angle complex $\zk$ is simply
  connected.
\end{corollary}
\begin{proof}
  Indeed, the 1-skeleton of our cellular decomposition of $\zk$ is contained
  in the torus $T^m=\rho^{-1}(1,\ldots,1)$.
\end{proof}

The cohomology of cellular pair $(\zk,T^m)$ also can be calculated by means of
the cochain complex $\mathcal C^{*}(K)$. The cellular cochain subcomplex
$\mathcal C^{*}(T^m)\subset\mathcal C^{*}(K)$ consists of monomials $u_I$
(i.e. monomials that do not contain $v_i$'s). This, of course, is just the
exterior algebra $\Lambda[u_1,\ldots,u_m]$. Hence,
\begin{equation}\label{pair}
  \mathcal C^{*}(\zk,T^m)=\mathcal C^{*}(\zk)/\Lambda[u_1,\ldots,u_m]
\end{equation}
(as complexes, not as algebras). We can also introduce relative bigraded
Betti numbers
\begin{equation}\label{rbbn}
  b_{-q,2p}(\zk,T^m)=\dim H^{-q,2p}[\mathcal C^{*}(\zk,T^m),d],
  \quad q,p=0,\ldots,m,
\end{equation}
define the $p$th relative Euler characteristic $\chi_p(\zk,T^m)$ as the Euler
number of complex $\mathcal C^{*,2p}(\zk,T^m)$:
\begin{equation}\label{rchip}
  \chi_p(\zk,T^m)=\sum_{q=0}^m(-1)^q\dim\mathcal C^{-q,2p}(\zk,T^m)
  =\sum_{q=0}^m(-1)^qb_{-q,2p}(\zk,T^m),
\end{equation}
and define the generating polynomial $\chi(\zk,T^m,t)$ as
$$
  \chi(\zk,T^m,t)=\sum_{p=0}^m\chi_p(\zk,T^m)t^{2p}.
$$
We will use the following theorem in the next section.

\begin{theorem}
\label{rgp}
  For any $(n-1)$-dimensional simplicial complex $K$ with $m$ vertices holds
  \begin{equation}
  \label{rhchi}
    \chi(\zk,T^m,t)=(1-t^2)^{m-n}(h_0+h_1t^2+\cdots+h_nt^{2n})-(1-t^2)^m.
  \end{equation}
\end{theorem}
\begin{proof}
Since $\mathcal C^{*}(T^m)=\Lambda[u_1,\ldots,u_m]$ and $\bideg u_i=(-1,2)$,
we have
$$
  \dim\mathcal C^{-q}(T^m)=\dim\mathcal C^{-q,2q}(T^m)=\textstyle\binom mq.
$$
It follows from (\ref{pair}), (\ref{chip}) and (\ref{rchip}) that
$$
  \chi_p(\zk,T^m)=\chi_p(\zk)-(-1)^p\dim\mathcal C^{-p,2p}(T^m).
$$
Hence,
\begin{align*}
  \chi(\zk,T^m,t)&=\chi(\zk,t)-\sum_{p=0}^m(-1)^p{\textstyle\binom
  mpt^{2p}}\\
  &=(1-t^2)^{m-n}(h_0+h_1t^2+\cdots+h_nt^{2n})-(1-t^2)^m,
\end{align*}
by (\ref{hchi}).
\end{proof}

At the end of this section we review the most important additional properties
of $\zk$ in the case when $|K|\cong S^{n-1}$, i.e. $K$ is a simplicial
sphere.

\begin{lemma}\label{manif}
  If $K$ is a simplicial sphere, i.e. $|K|=S^{n-1}$, then $\zk$ is an
  $(m+n)$-dimensional (closed) manifold.\ep
\end{lemma}
In~\cite[p. 434]{DJ} the authors considered the manifold $\mathcal Z$ defined
for any simple $n$-polytope $P^{*}$ with $m$ facets as $\mathcal Z=(T^m\times
P^{*})/\,\sim$, where $\sim$ is a certain equivalence relation. We showed
in~\cite{BP2} that if $K$ is a polytopal sphere, i.e. $K=\partial P^n$ for
some simplicial polytope $P^n$, then our moment-angle complex $\zk$ coincides
with the manifold $\mathcal Z$ defined by simple polytope $P^{*}$ dual to
$P$. For the case of general simplicial sphere $K$,
see~\cite{BP3},~\cite{BP4}.

\begin{theorem}\label{fc}
  Let $K$ be an $(n-1)$-dimensional simplicial sphere, and let $\zk$ be the
  corresponding moment-angle manifold. The fundamental cohomological class of
  $\zk$ is represented by any monomial $\pm v_Iu_J\in\mathcal C(K)$ of
  bidegree $(-(m-n),2m)$ such that $I$ is an $(n-1)$-simplex of $K$ and
  $I\cap J=\varnothing$. The choice of sign depends on the orientation of
  $\zk$.
\end{theorem}
\begin{proof}
We have $\dim\zk=m+n$. It follows from Lemma~\ref{bbgen}~(f) that
$H^{m+n}(\zk)=H^{-(m-n),2m}(\zk)$. By definition, the module $\mathcal
C^{-(m-n),2m}(\zk)$ is spanned by monomials $v_Iu_J$ such that $I\in
K^{n-1}$, $\#I=n$, $J=[m]\setminus I$, and all these monomials are cocycles.
Suppose that $I,I'$ are two $(n-1)$-simplices of $K^{n-1}$ sharing a common
$(n-2)$-face. Then the corresponding cocycles $v_Iu_J$, $v_{I'}u_{J'}$, where
$J=[m]\setminus I$, $J'=[m]\setminus I'$, are cohomological (up to sign).
Indeed, let $v_Iu_J=v_{i_1}\cdots v_{i_n}u_{j_1}\cdots u_{j_{m-n}}$,
$v_{I'}u_{J'}= v_{i_1}\cdots v_{i_{n-1}}v_{j_1}u_{i_n}u_{j_2}\cdots
u_{j_{m-n}}$. Since any $(n-2)$-face of $K$ is contained in exactly two
$(n-1)$-faces, the identity
\begin{multline*}
  d(v_{i_1}\cdots v_{i_{n-1}}u_{i_n}u_{j_1}u_{j_2}\cdots u_{j_{m-n}})\\
  =v_{i_1}\cdots v_{i_n}u_{j_1}\cdots u_{j_{m-n}}-
  v_{i_1}\cdots v_{i_{n-1}}v_{j_1}u_{i_n}u_{j_2}\cdots u_{j_{m-n}}
\end{multline*}
holds in $\mathcal C(K)\subset\k(K)\otimes\Lambda[u_1,\ldots,u_m]$, hence,
$v_Iu_J$ and $v_{I'}u_{J'}$ are cohomological. Since $K^{n-1}$ is a
simplicial sphere, any two $(n-1)$-simplices of $K^{n-1}$ can be connected by
a chain of simplices such that any two successive simplices share a common
$(n-2)$ face. Thus, any two cocycles in $\mathcal C^{-(m-n),2m}(\zk)$ are
cohomological, and we can take any one as a representative for the
fundamental cohomological class of $\zk$ (after a proper choice of sign).
\end{proof}

\begin{remark}
In the proof of the above theorem we have used two combinatorial properties
of $K^{n-1}$. The first one is that any $(n-2)$-face is contained in exactly
two $(n-1)$-faces, and the second one is that any two
$(n-1)$-simplices can be connected by a chain of simplices such
that any two successive simplices share a common $(n-2)$-face. Both
properties hold for any simplicial manifold. Hence, for any
simplicial manifold $K^{n-1}$ one has $b_{m+n}(\zk)=b_{-(m-n),2m}(\zk)=1$ and
the generator of $H^{m+n}(\zk)$ can be chosen as in Theorem~\ref{fc}.
\end{remark}

\begin{corollary}\label{bpd}
  The Poincar\'e duality for the moment angle manifold $\zk$ defined by a
  simplicial sphere $K^{n-1}$ regards the bigraded structure in the
  (co)homology, i.e.
  $$
    H^{-q,2p}(\zk)\cong H_{-(m-n)+q,2(m-p)}(\zk).
  $$
  In particular,
  \begin{equation}\label{bpdbn}
    b_{-q,2p}(\zk)=b_{-(m-n)+q,2(m-p)}(\zk).\qquad\square
  \end{equation}
\end{corollary}
\begin{corollary}\label{bbss}
  Let $K^{n-1}$ be an $(n-1)$-dimensional simplicial sphere, and let $\zk$
  be the corresponding moment-angle complex, $\dim\zk=m+n$. Then
  \begin{itemize}
  \item[\it (a)] $b_{-q,2p}(\zk)=0$ if $q\ge m-n$, with only exception
   $b_{-(m-n),2m}=1$;
  \item[\it (b)] $b_{-q,2p}(\zk)=0$ if $p-q\ge n$, with only exception
   $b_{-(m-n),2m}=1$.\ep
  \end{itemize}
\end{corollary}
It follows that if $K^{n-1}$ is a simplicial sphere, then non-zero
bigraded Betti numbers $b_{r,2p}(\zk)$, $r\ne0$, $r\ne m-n$, appear only in
the ``strip" bounded by the lines $r=-(m-n-1)$, $r=-1$, $p+r=1$ and $p+r=n-1$
in the second quadrant (see Figure~4~(b)). Compare this with Figure~4~(a)
corresponding to the case of general $K$.

It follows from (\ref{chip}) and (\ref{bpdbn}) that for any simplicial sphere
$K$ holds
$$
  \chi_p(\zk)=(-1)^{m-n}\chi_{m-p}(\zk).
$$
From this and (\ref{hchi}) we get
\begin{multline*}
  \frac{h_0+h_1t^2+\cdots+h_nt^{2n}}{(1-t^2)^n}=
  (-1)^{m-n}\frac{\chi_m+\chi_{m-1}t^2+\cdots+\chi_0t^{2m}}{(1-t^2)^m}\\=
  (-1)^n\frac{\chi_0+\chi_1t^{-2}+\cdots+\chi_mt^{-2m}}{(1-t^{-2})^m}
  =(-1)^n\frac{h_0+h_1t^{-2}+\cdots+h_nt^{-2n}}{(1-t^{-2})^n}\\=
  \frac{h_0t^{2n}+h_1t^{2(n-1)}+\cdots+h_n}{(1-t^2)^n}.
\end{multline*}
Hence, $h_i=h_{n-i}$. Thus, we have deduced the Dehn--Sommerville equations
as a corollary of the bigraded Poincar\'e duality~(\ref{bpdbn}).

The identity~(\ref{hchi}) also allows to interpret different inequalities for
the numbers of faces of simplicial spheres (resp. simplicial manifolds) in
terms of topological invariants (bigraded Betti numbers) of the corresponding
moment-angle manifolds (resp. complexes) $\zk$.

\begin{example}
It follows from Lemma~\ref{bbgen} that for any $K$ holds
\begin{align*}
  &\chi_0(K)=1,\\
  &\chi_1(K)=0,\\
  &\chi_2(K)=-b_{-1,4}(\zk)=-b_3(\zk),\\
  &\chi_3(K)=b_{-2,6}(\zk)-b_{-1,6}(\zk)
\end{align*}
(note that $b_4(\zk)=b_{-2,6}(\zk)$, while
$b_5(\zk)=b_{-1,6}(\zk)+b_{-3,8}(\zk)$). Now, identity~(\ref{hchi}) shows
that
\begin{align*}
  &h_0=1,\\
  &h_1=m-n,\\
  &h_2={\textstyle\binom{m-n+1}2}-b_3(\zk),\\
  &h_3={\textstyle\binom{m-n+2}3}-(m-n)b_{-1,4}(\zk)+
   b_{-2,6}(\zk)-b_{-1,6}(\zk).
\end{align*}
It follows that the inequality $h_1\le h_2$, $n\ge4$, from the Generalized
Lower Bound hypothesis~(\ref{glb}) for simplicial spheres is equivalent to
the following inequality:
\begin{equation}
\label{h12}
  b_3(\zk)\le\textstyle\binom{m-n}2.
\end{equation}
The next inequality $h_2\le h_3$, $n\ge6$, from~(\ref{glb}) is equivalent to
the following inequality for the bigraded Betti numbers of $\zk$:
\begin{equation}
\label{h23}
  {\textstyle\binom{m-n+1}3}-(m-n-1)b_{-1,4}(\zk)+
  b_{-2,6}(\zk)-b_{-1,6}(\zk)\ge0.
\end{equation}
\end{example}

Thus, we see that the combinatorial Generalized Lower Bound inequalities are
interpreted as ``topological" inequalities for the (bigraded) Betti numbers
of a certain manifold. So, one can try to use topological methods (such as
the equivariant topology or Morse theory) to prove inequalities
like~(\ref{h12}) or~(\ref{h23}).  Such topological approach to the hypotheses
like $g$-theorem or Generalized Lower Bound has an advantage of being
independent on whether the simplicial sphere $K$ is polytopal or not. Indeed,
all known proofs of the necessity of $g$-theorem for simplicial polytopes
(including the original one by Stanley~\cite{St2}, McMullen's
proof~\cite{McM}, and the recent proof by Timorin~\cite{Ti}) follow the same
scheme. Namely, the numbers $h_i$, $i=1,\ldots,n$, are interpreted as the
dimensions of graded components $A^i$ of a certain graded algebra $A$
satisfying the Hard Lefschetz Theorem. The latter means that there is an
element $\omega\in A^1$ such that the multiplication by $\omega$ defines a
{\it monomorphism} $A^i\to A^{i+1}$ for $i<\bigl[\frac n2\bigr]$. This
implies $h_i\le h_{i+1}$ for $i<\bigl[\frac n2\bigr]$. However, such element
$\omega$ is lacking for non-polytopal $K$, which means that one should
develop a new technique for proving the $g$-theorem (or, better to say,
$g$-conjecture) for simplicial spheres. Certainly, it may happen that the
$g$-theorem fails to be true for simplicial spheres, however, many recent
efforts of computer-aided seek for counter examples were unsuccessful (see,
e.g.,~\cite{BjLu}).

It can be easily seen that our moment-angle complex $\zk$ is a manifold
provided that the cone $\cone(K)$ is non-singular. This is equivalent to the
condition that the suspension $\Sigma|K|$ is a (topological) manifold, which
implies (due to the suspension isomorphism and Poincar\'e duality in the
homology) that $|K|$ is a homology sphere. An important class of simplicial
homology spheres is known in combinatorics as Gorenstein* complexes (see,
e.g.,~\cite{St3} for definition).  As it was pointed out by Stanley
in~\cite{St3}, the Gorenstein* complexes are the most general objects
appropriate for generalizing the $g$-theorem (they include polytopal spheres,
PL-spheres and simplicial spheres as particular cases). In our terms, the
Gorenstein* complexes $K$ (see~\cite[p.75]{St1}) can be characterized by the
condition that the $\Tor$-algebra
$\Tor_{\k[v_1,\ldots,v_m]}\bigl(\k(K),\k\bigr)$ satisfies the bigraded
duality~(\ref{bpdbn}), i.e. $\zk$ is a Poincar\'e duality space (not
necessarily a manifold). In particular, the Dehn--Sommerville relations
$h_i=h_{n-i}$ continue to hold for Gorenstein* complexes.

\section{Homology of $\mathcal W_K$ and generalized Dehn--Sommerville
equations}

Here we assume that $K^{n-1}$ is a triangulation of a manifold, i.e., a
simplicial manifold. In this case the moment-angle complex $\zk$ is not a
manifold, however, its singularities can be easily treated. Indeed,
$|\cc(K)|$ is homeomorphic to $|\cone(K)|$, and the vertex of the cone is the
point $p=(1,\ldots,1)\in|\cc(K)|\subset I^m$.  Let
$U_\varepsilon(p)\subset|\cc(K)|$ be a small neighbourhood of $p$ in
$|\cc(K)|$. Then the closure of $U_\varepsilon(p)$ is also homeomorphic to
$|\cone(K)|$. It follows from the definition of $\zk$ (see~(\ref{zkwk})) that
$U_\varepsilon(T^m)=\rho^{-1}\bigl(U_\varepsilon(p)\bigr)\subset\zk$ is a
small invariant neighbourhood of the torus $T^m=\rho^{-1}(p)$ in $\zk$. Here
$\rho:(D^2)^m\to I^m$ is the orbit map. Then for small $\varepsilon$ the
closure of the neighbourhood $U_\varepsilon(T^m)$ is homeomorphic to
$|\cone(K)|\times T^m$. Taking $U_\varepsilon(T^m)$ away from $\zk$ we obtain
a manifold with boundary, which we denote $W_K$. Hence,
we have
$$
  W_K=\overline{\zk\setminus|\cone(K)|\times T^m},\quad
  \partial W_K=|K|\times T^m.
$$
Note that since the neighbourhood $U_\varepsilon(T^m)$ is $T^m$-invariant,
the torus $T^m$ acts on $W_K$.

\begin{theorem}
\label{homotwk}
  The manifold with boundary $W_K$ is equivariantly homotopy equivalent to
  the moment-angle complex $\wk$ (see~{\rm(\ref{zkwk})}). There is a
  canonical relative isomorphism of pairs $(W_K,\partial W_K)\to(\zk,T^m)$.
\end{theorem}
\begin{proof}
To prove the first assertion we construct homotopy equivalence
$|\cc(K)|\setminus U_\varepsilon(p)\to|\cub(K)|$ (see Proposition~\ref{cubk})
as it is shown on Figure~5. This homotopy equivalence is covered by a
$T^m$-invariant homotopy equivalence $W_K=\zk\setminus
U_\varepsilon(T^m)\to\wk$, as needed. The second assertion follows easily
from the definition of $W_K$.
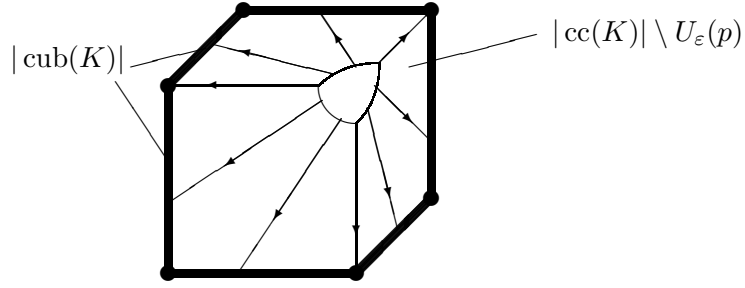
\begin{figure}
  \begin{picture}(120,45)
  \put(70,5){\circle*{2}}
  \put(80,15){\circle*{2}}
  \put(45,30){\circle*{2}}
  \put(55,40){\circle*{2}}
  \put(45,5){\circle*{2}}
  \put(80,40){\circle*{2}}
  \multiput(45,29.3)(0,0.1){16}{\line(1,1){10}}
  \multiput(70,4.3)(0,0.1){16}{\line(1,1){10}}
  \put(73,33){\vector(1,1){4}}
  \put(73,33){\line(1,1){7}}
  \put(45,30){\line(1,0){20}}
  \put(70,5){\line(0,1){20}}
  \put(65,30){\vector(-1,0){15}}
  \put(70,25){\vector(0,-1){15}}
  \put(65.5,28){\line(-3,-2){20.5}}
  \put(68,25.6){\line(-2,-3){13.5}}
  \put(65.5,28){\vector(-3,-2){13}}
  \put(68,25.6){\vector(-2,-3){9}}
  \put(67,31.5){\line(-4,1){16}}
  \put(70,32.5){\line(-2,3){5}}
  \put(67,31.5){\vector(-4,1){12.5}}
  \put(70,32.5){\vector(-2,3){3}}
  \put(71.5,27){\line(1,-4){4}}
  \put(72.5,30){\line(1,-1){7}}
  \put(71.5,27){\vector(1,-4){3}}
  \put(72.5,30){\vector(1,-1){5}}
  \put(70,30){\oval(10,10)[lb]}
  \qbezier(65,30)(68,33)(73,33)
  \qbezier(70,25)(73,28)(73,33)
  \put(77,32.5){\line(4,1){17}}
  \put(96,36){$|\cc(K)|\setminus U_\varepsilon(p)$}
  \put(45,20){\line(-2,3){7}}
  \put(50,35){\line(-4,-1){10}}
  \put(24,32.5){$|\cub(K)|$}
  \linethickness{1mm}
  \put(45,5){\line(1,0){25}}
  \put(55,40){\line(1,0){25}}
  \put(45,5){\line(0,1){25}}
  \put(80,15){\line(0,1){25}}
  \end{picture}
  \caption{Homotopy equivalence
  $|\cc(K)|\setminus U_\varepsilon(p)\to|\cub(K)|$.}
  \end{figure}
\end{proof}

According to Lemma \ref{macell}, the moment-angle complex $\wk\subset(D^2)^m$
has a cellular structure with 5 different cell types $D_i$, $I_i$, $0_i$,
$T_i$, $1_i$, $i=1,\ldots,m$, (see Figure~3). The homology of $\wk$ (and of
$W_K$) can be calculated by means of the corresponding cellular chain
complex, which we denote $\bigr[\mathcal C_{*}(\wk),\partial_c\bigl]$. In
comparison with the moment-angle complex $\zk$ studied in the previous
section the complex $\wk$ has more types of cells (remember that $\zk$ has
only 3 cell types $D_i$, $T_i$, $1_i$). However, the wonderful thing is that
the cellular chain complex $\bigr[\mathcal C_{*}(\wk),\partial_c\bigl]$ can
be canonically made {\it bigraded}. Namely, the following statement holds
(compare with~(\ref{bgcellz})).
\begin{lemma}
\label{bgw}
  Put
  \begin{gather}
  \label{bgcellw}
    \bideg{D_i}=(0,2),\quad\bideg{T_i}=(-1,2),\quad\bideg{I_i}=(1,0),\\
    \bideg{0_i}=\bideg{1_i}=(0,0),\quad i=1,\ldots,m.\notag
  \end{gather}
  This makes the cellular chain complex $\bigr[\mathcal
  C_{*}(\wk),\partial_c\bigl]$ a bigraded differential module with
  differential $\partial_c$ adding $(-1,0)$ to bidegree. The original
  grading of $\mathcal C_{*}(\wk)$ by dimensions of cells corresponds to the
  total degree (i.e., the dimension of a cell equals the sum of its two
  degrees).
\end{lemma}
\begin{proof}
We need only to check that the differential $\partial_c$ adds $(-1,0)$ to
bidegree. This follows from~(\ref{bgcellw}) and
$$
  \partial_cD_i=T_i,\quad\partial_cI_i=1_i-0_i,\quad
  \partial_cT_i=\partial_c1_i=\partial_c0_i=0.
$$
\end{proof}
Note that, unlike bigraded structure in $\mathcal C_{*}(\zk)$, elements of
$\mathcal C_{*,*}(\wk)$ may have {\it positive} first degree (due to the
positive first degree of $I_i$'s). However, as in the case of $\zk$, the
differential $\partial_c$ does not change the second degree, which allows to
split the bigraded complex $\mathcal C_{*,*}(\wk)$ to the sum of complexes
$\mathcal C_{*,2p}(\wk)$, $p=0,\ldots,m$.

In the same way as we have done this for $\zk$ and for the pair $(\zk,T^m)$
define the bigraded Betti numbers
\begin{equation}\label{bbnw}
  b_{q,2p}(\wk)=\dim H_{q,2p}[\mathcal C_{*,*}(\wk),\partial_c],
  \quad -m\le q\le m,\;0\le p\le m
\end{equation}
(note that $q$ may be both positive and negative), the $p$th Euler
characteristic $\chi_p(\wk)$ as the Euler number of complex $\mathcal
C_{*,2p}(\wk)$:
\begin{equation}\label{chipw}
  \chi_p(\wk)=\sum_{q=-m}^m(-1)^q\dim\mathcal C_{q,2p}(\wk)
  =\sum_{q=-m}^m(-1)^qb_{q,2p}(\wk),
\end{equation}
and the generating polynomial $\chi(\wk,t)$ as
$$
  \chi(\wk,t)=\sum_{p=0}^m\chi_p(\wk)t^{2p}.
$$

The following theorem provides the exact formula for generating polynomial
$\chi(\wk,t)$ and is analogous to theorems~\ref{gpz} and~\ref{rgp}

\begin{theorem}
\label{chitwk}
  For any $(n-1)$-dimensional simplicial complex $K$ with $m$ vertices holds
  \begin{align*}
    \chi(\wk,t)&=
    (1-t^2)^{m-n}(h_0+h_1t^2+\cdots+h_nt^{2n})+
    \bigl(\chi(K)-1\bigr)(1-t^2)^m\\
    &=(1-t^2)^{m-n}(h_0+h_1t^2+\cdots+h_nt^{2n})+(-1)^{n-1}h_n(1-t^2)^m,
  \end{align*}
  where $\chi(K)=f_0-f_1+\ldots+(-1)^{n-1}f_{n-1}=1+(-1)^{n-1}h_n$ is the
  Euler number of $K$.
\end{theorem}
\begin{proof}
Lemma \ref{macell} shows that every moment-angle complex $\ma(C)$ is a
cellular subcomplex of $(D^2)^m$, and each cell is the product of cells of 5
different types: $D_i$, $I_i$, $0_i$, $T_i$, and $1_i$, $i=1,\ldots,m$. These
products are encoded by words $D_II_J0_LT_P1_Q$, where $I,J,L,P,Q$ are
pairwise disjoint subsets of $[m]$ such that $I\cup J\cup K\cup P\cup Q=[m]$.
In the case $\ma(C)=\wk=\ma\bigl(\cub(K)\bigr)$ the definition of $\cub(K)$
(see Proposition~\ref{cubk}) shows that the cell $D_II_J0_LT_P1_Q$ belongs to
$\wk$ if and only if the following two conditions are satisfied:
\begin{enumerate}
\item The set $I\cup J\cup L$ is a simplex of $K^{n-1}$.
\item $\#L\ge1$.
\end{enumerate}
Let $c_{ijlpq}(\wk)$ denote the number of cells $D_II_J0_LT_P1_Q\subset\wk$
with $i=\#I$, $j=\#J$, $l=\#L$, $p=\#P$, $q=\#Q$, $i+j+l+p+q=m$. It follows
that
\begin{equation}
\label{ci}
  c_{ijlpq}(\wk)=f_{i+j+l-1}\textstyle\binom{i+j+l}i\binom{j+l}l
  \binom{m-i-j-l}p,
\end{equation}
where $(f_0,\ldots,f_{n-1})$ is the $f$-vector of $K$ (we also assume
$f_{-1}=1$ and $f_k=0$ for $k<-1$ or $k>n-1$). By (\ref{bgcellw}),
$$
  \bideg(D_II_J0_LT_P1_Q)=(j-p,2(i+p)).
$$
Now we calculate $\chi_r(\wk)$ as it is defined by (\ref{chipw}), using
(\ref{ci}):
$$
  \chi_r(\wk)=\mathop{\sum_{i,j,l,p}}\limits_{i+p=r,l\ge1}(-1)^{j-p}
  f_{i+j+l-1}\textstyle\binom{i+j+l}i\binom{j+l}l\binom{m-i-j-l}p.
$$
Substituting $s=i+j+l$, $i=r-p$ above, we obtain
\begin{align*}
  \chi_r(\wk)&=\mathop{\sum_{l,s,p}}\limits_{l\ge1}(-1)^{s-r-l}
  f_{s-1}\textstyle{\binom{s}{r-p}\binom{s-r+p}l\binom{m-s}p}\\
  &=\sum_{s,p}\Bigl((-1)^{s-r}f_{s-1}{\textstyle\binom{s}{r-p}\binom{m-s}p}
  \sum_{l\ge1}(-1)^l{\textstyle\binom{s-r+p}l}\Bigl)
\end{align*}
Since
$$
  \sum_{l\ge1}(-1)^l{\textstyle\binom{s-r+p}l}=
  \left\{
  \begin{aligned}
    -1,&\quad s>r-p,\\
    0,&\quad s\le r-p,
  \end{aligned}
  \right.
$$
we obtain
\begin{align*}
  \chi_r(\wk)&=-\mathop{\sum_{s,p}}\limits_{s>r-p}(-1)^{s-r}
  f_{s-1}{\textstyle\binom{s}{r-p}\binom{m-s}p}\\
  &=-\sum_{s,p}(-1)^{r-s}
  f_{s-1}{\textstyle\binom{s}{r-p}\binom{m-s}p}+
  \sum_s(-1)^{r-s}f_{s-1}{\textstyle\binom{m-s}{r-s}}.
\end{align*}
The second sum in the above formula is exactly $\chi_r(\zk)$
(see~(\ref{chipzk})). To calculate the first sum, we mention that
$$
  \sum_p\textstyle\binom{s}{r-p}\binom{m-s}p=\binom mr
$$
(this follows from calculating the coefficient of $\alpha^r$ in both sides of
$(1+\alpha)^s(1+\alpha)^{m-s}=(1+\alpha)^m$). Hence,
$$
  \chi_r(\wk)=-\sum_s(-1)^{r-s}f_{s-1}{\textstyle\binom mr}+\chi_r(\zk)=
  (-1)^r{\textstyle\binom mr}\bigl( \chi(K)-1 \bigr)+\chi_r(\zk),
$$
since $-\sum_s(-1)^sf_{s-1}=\chi(K)-1$ (remember that $f_{-1}=1$). Finally,
using~(\ref{hchi}), we calculate
\begin{align*}
  \chi(\wk,t)=\sum_{r=0}^m\chi_r(\wk)t^{2r}=
  \sum_{r=0}^m(-1)^r{\textstyle\binom mr}\bigl( \chi(K)-1 \bigr)t^{2r}+
  \sum_{r=0}^m\chi_r(\zk)t^{2r}\\
  =\bigl( \chi(K)-1 \bigr)(1-t^2)^m+
  (1-t^2)^{m-n}(h_0+h_1t^2+\cdots+h_nt^{2n}).
\end{align*}
\end{proof}

Suppose now that $K$ is an orientable simplicial manifold. It is easy
to see that then $W_K$ is also orientable. Hence, there
are relative Poincar\'e duality isomorphisms:
\begin{equation}\label{rpd}
  H_k(W_K)\cong H^{m+n-k}(W_K,\partial_c W_K), \quad k=0,\ldots,m.
\end{equation}

\begin{corollary}
\label{DSsm}
The following relations hold for the $h$-vector $(h_0,h_1,\ldots,h_n)$ of any
simplicial manifold $K^{n-1}$:
$$
  h_{n-i}-h_i=(-1)^i\bigl(\chi(K^{n-1})-\chi(S^{n-1})\bigr)
  {\textstyle\binom ni},\quad i=0,1,\ldots,n,
$$
where $\chi(S^{n-1})=1+(-1)^{n-1}$ is the Euler number of an $(n-1)$-sphere.
\end{corollary}
\begin{proof}
Theorem~\ref{homotwk} shows that $H^{m+n-k}(W_K,\partial_c W_K)=
H^{m+n-k}(\zk,T^m)$ and $H_k(W_K)=H_k(\wk)$.  Moreover, it can be seen in the
same way as in Corollary~\ref{bpd} that relative Poincar\'e duality
isomorphisms~(\ref{rpd}) regard the bigraded structures in the (co)homology
of $\wk$ and $(\zk,T^m)$. It follows that
$$
  b_{-q,2p}(\wk)=b_{-(m-n)+q,2(m-p)}(\zk,T^m).
$$
Hence,
$$
  \chi_p(\wk)=(-1)^{m-n}\chi_{m-p}(\zk,T^m),
$$
and
\begin{equation}
\label{relchid}
  \begin{aligned}
    \chi(\wk,t)&=(-1)^{m-n}\sum_p\chi_{m-p}(\zk,T^m)t^{2p}\\
    &=(-1)^{m-n}t^{2m}\chi(\zk,T^m,\textstyle\frac1t).
  \end{aligned}
\end{equation}
Using (\ref{rhchi}), we calculate
\begin{multline*}
  (-1)^{m-n}t^{2m}\chi(\zk,T^m,{\textstyle\frac1t})\\
      =(-1)^{m-n}t^{2m}(1-t^{-2})^{m-n}(h_0+h_1t^{-2}+\cdots+h_nt^{-2n})\\
      -(-1)^{m-n}t^{2m}(1-t^{-2})^m\\
    =(1-t^2)^{m-n}(h_0t^{2n}+h_1t^{2n-2}+\cdots+h_n)+
    (-1)^{n-1}(1-t^2)^m.
\end{multline*}
Substituting the formula for $\chi(\wk,t)$ from Theorem~\ref{chitwk} and the
above expression into formula~(\ref{relchid}), we obtain
\begin{multline*}
  (1-t^2)^{m-n}(h_0+h_1t^2+\cdots+h_nt^{2n})+\bigl(\chi(K)-1\bigr)
  (1-t^2)^m\\=(1-t^2)^{m-n}(h_0t^{2n}+h_1t^{2n-2}+\cdots+h_n)+
  (-1)^{n-1}(1-t^2)^m.
\end{multline*}
Calculating the coefficient of $t^{2i}$ in both sides after dividing the
above identity by $(1-t^2)^{m-n}$, we obtain
$h_{n-i}-h_i=(-1)^i\bigl(\chi(K^{n-1})-\chi(S^{n-1})\bigr)
{\textstyle\binom ni}$, as needed.
\end{proof}
Corollary~\ref{DSsm} generalize the Dehn--Sommerville
equations~(\ref{DSpol}) for simplicial spheres. If $|K|=S^{n-1}$ or $(n-1)$
is odd, Corollary~\ref{DSsm} gives just $h_{n-i}=h_i$.

\begin{corollary}
The following relations hold for any simplicial manifold~$K^{n-1}$:
$$
  h_{n-i}-h_i=(-1)^i(h_n-1){\textstyle\binom ni},\quad i=0,1,\ldots,n.
$$
\end{corollary}
\begin{proof}
Since
$\chi(K^{n-1})=1+(-1)^{n-1}h_n$, $\chi(S^{n-1})=1+(-1)^{n-1}$, we have
$$
  \chi(K^{n-1})-\chi(S^{n-1})=(-1)^{n-1}(h_n-1)=(h_n-1)
$$
(the coefficient $(-1)^{n-1}$ can be dropped since for odd $(n-1)$
the left side is zero).
\end{proof}

\begin{corollary}
  For any $(n-1)$-dimensional simplicial manifold the numbers $h_{n-i}-h_i$,
  $i=0,1,\ldots,n$, are homotopy invariants. In particular, they are
  independent on a triangulation.
\end{corollary}
In the particular case of PL-manifolds the topological invariance of numbers
$h_{n-i}-h_i$ was firstly observed by Pachner in~\cite[(7.11)]{Pa}.

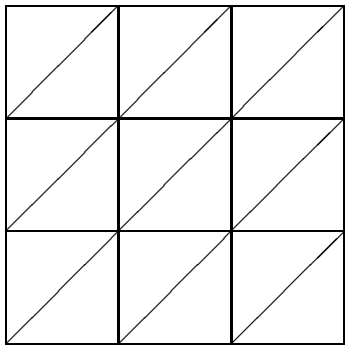
\begin{figure}
  \begin{picture}(120,45)(-40,0)
  \multiput(0,0)(0,15){4}{\line(1,0){45}}
  \multiput(0,0)(15,0){4}{\line(0,1){45}}
  \put(0,30){\line(1,1){15}}
  \put(0,15){\line(1,1){30}}
  \put(0,0){\line(1,1){45}}
  \put(15,0){\line(1,1){30}}
  \put(30,0){\line(1,1){15}}
  \end{picture}
  \caption{Triangulation of $T^2$ with $f=(9,27,18)$, $h=(1,6,12,-1)$.}
\end{figure}
\begin{example}
Consider triangulations of the 2-torus $T^2$, so $n=3$, $\chi(T^2)=0$. Since
for any $K^{n-1}$ holds $\chi(K^{n-1})=1+(-1)^{n-1}h_n$, in our case we
have $h_3=-1$. Then Corollary~\ref{DSsm} gives
$$
  h_3-h_0=-2,\quad h_2-h_1=6.
$$
For instance, the triangulation on Figure~6 has $f_0=9$ vertices, $f_1=27$
edges and $f_2=18$ triangles. The corresponding $h$-vector is $(1,6,12,-1)$.
\end{example}

\section{Concluding remarks}

The main goal of the present paper was to establish new connections between
topology of manifolds and cellular complexes and combinatorics by means of
the notion of a moment-angle complex, introduced by the authors in previous
papers~\cite{BP2},~\cite{BP3},~\cite{BP4}. As we have seen, the combinatorics
of simplicial manifolds and related objects (polytopes, simplicial spheres,
simplicial complexes, coordinate subspace arrangements) can be effectively
described by means of topological invariants of bigraded equivariant
moment-angle complexes. One of the main properties of a moment-angle complex
is the existence of a torus action all of whose isotropy subgroups are
coordinate ones. This, in particular, allows to introduce an additional
grading to the cohomology ring of the moment-angle complex. On this point one
can observe that there is a natural $\Z/2$-analogue of almost all
constructions presented in this paper. The first step is to replace the torus
$T^m$ by its ``real analogue", namely, the group $(\Z/2)^m$. Then the unit
cube $I^m=[0,1]^m$ is the orbit space for the action of $(\Z/2)^m$ on the
bigger cube $[-1,1]^m$, the ``real analogue" of the poly-disk
$(D^2)^m\subset\C^m$. Now, starting from any cubical subcomplex $C\subset
I^m$ one can construct another $(\Z/2)^m$-symmetrical cubical complex
embedded into $[-1,1]^m\subset\R^m$, just in the same way as we did it in
Definition~\ref{ma}. In particular, for any simplicial complex $K$ on $m$
vertices one can construct $(\Z/2)^m$-symmetrical cubical complexes $\zk^\R$,
$\wk^\R$, the ``real analogues" of moment-angle complexes $\zk$, $\wk$,
see~(\ref{zkwk}). The analogue of Lemma~\ref{he} holds for {\it real}
coordinate subspace arrangements. Namely, the complement $U\bigl(\mathcal
A^\R(K)\bigr)$ of the real coordinate subspace arrangement $\mathcal A^\R(K)$
defined by $K$ in the same way as in~(\ref{compl}) is
$(\Z/2)^m$-equivariantly homotopy equivalent to $\zk^\R$. However, the
situation with the cohomology algebra of $\zk^\R$ is more subtle: as we have
already mentioned, the analogue of Theorem~\ref{arrang} does not hold for
$\zk^\R$, at least for $\Q$-coefficients (this is a usual thing in topology:
the cohomology of ``real" objects is more complicated than that of ``complex"
ones). At the same time $\zk^\R$ is an $m$-dimensional manifold provided that
$K$ is a simplicial sphere. So, for any simplicial sphere $K$ with $m$
vertices we have a $(\Z/2)^m$-symmetric manifold with $(\Z/2)^m$-invariant
cubical complex structure. This class of cubical manifolds may be useful in
the cubical analogue of the combinatorial theory of $f$-vectors of simplicial
complexes.

\end{document}